\documentclass[AMA,STIX1COL]{WileyNJD-v2}
\usepackage{booktabs}
\usepackage{textcomp}
\usepackage{multirow}
\usepackage{verbatim}
\usepackage{amsfonts}
\usepackage{amsmath}
\usepackage{cancel}
\usepackage{mathrsfs}
\usepackage{graphicx}
\usepackage{tikz}
\usepackage{pgfplots}
\usepackage{float}
\usetikzlibrary{pgfplots.groupplots}
\usetikzlibrary{decorations.pathreplacing,calligraphy}
\usepackage{mathtools}
\usepackage{geometry}
\usepackage{url}
\usepackage{siunitx}
\usepackage{tabularx}
\usepackage{caption}
\usepackage{bm}
\usepackage{moreverb}

\newcommand{\dpart}[2]{\displaystyle\frac{\partial #1}{\partial #2}}

\makeatother

\received{}
\revised{}
\accepted{}

\newcommand\BibTeX{{\rmfamily B\kern-.05em \textsc{i\kern-.025em b}\kern-.08em
T\kern-.1667em\lower.7ex\hbox{E}\kern-.125emX}}

\articletype{Research Article}%
\pgfplotsset{compat=1.18}

\begin{document}

\title{Efficient hyperbolic-parabolic models on multi-dimensional unbounded domains using an extended DG approach}

\author[CASA]{Federico Vismara}
\author[LDO,DMI]{Tommaso Benacchio}

\authormark{Federico Vismara and Tommaso Benacchio} 

\address[CASA]{CASA, Technische Universiteit Eindhoven, The Netherlands}
\address[LDO]{Leonardo Labs, Future Rotorcraft Technologies, Leonardo S.p.A., Italy}
\address[DMI]{Now at: Danish Meteorological Institute, Lyngbyvej 100, 2100 Copenhagen, Denmark}

\corres{Federico Vismara, CASA, Technische Universiteit Eindhoven, The Netherlands. \email{f.vismara@tue.nl}}

\abstract{We introduce an extended discontinuous Galerkin discretization of hyperbolic-parabolic problems on multidimensional semi-infinite domains. Building on previous work on the one-dimensional case, we split the strip-shaped computational domain into a bounded region, discretized by means of discontinuous finite elements using Legendre basis functions, and an unbounded subdomain, where scaled Laguerre functions are used as a basis. Numerical fluxes at the interface allow for a seamless coupling of the two regions. The resulting coupling strategy is shown to produce accurate numerical solutions in tests on both linear and non-linear scalar and vectorial model problems. In addition, an efficient absorbing layer can be simulated in the semi-infinite part of the domain in order to damp outgoing signals with negligible spurious reflections at the interface. By tuning the scaling parameter of the Laguerre basis functions, the extended DG scheme  simulates transient dynamics over large spatial scales with a substantial reduction in computational cost at a given accuracy level compared to standard single-domain discontinuous finite element techniques.}


\maketitle

\section{Introduction}
\label{sec:intro} \indent

Accurate computational representation of physical phenomena on unbounded domains remains a challenge for numerical methods for transient differential problems, with fields of interest ranging from cell growth modelling \cite{xia:2020} to upper atmosphere dynamics \cite{klemp:2022} and space weather \cite{akmaev:2011,jackson:2019}. In the quest to make the inherently infinite problem tractable, numerical treatments based on absorbing (also called \emph{sponge}) layers are usually preferred to analytical approaches relying on non-reflecting boundary conditions \cite{dea:2011,engquist:1977,israeli:1981,neta:2008} (though see also \cite{dolci:2022} for a recent study challenging that convention). The former class of methods typically relies on dividing the computational domain into a bounded region of interest to the physical phenomena relevant to the problem at hand (e.g., the troposphere and lower stratosphere in currently operational atmospheric models \cite{melvin:2019,wood:2014}) and an unbounded buffer region for damping the outgoing signals. As the interface between the two subdomains is essentially artificial, the numerical treatment of the interface and the buffer region need to cause reflections in the bounded region that are as low-amplitude as possible. An optimal setup with absorbing boundary conditions maximizes the extent of the region of interest while keeping enough computational points in the buffer zone so as to minimize reflections due to the artificial boundary and the buffer zone itself. To that end, fine-tuning of the sponge layer parameters is usually required, making it challenging to obtain efficient numerical tools agnostic to the properties of the outgoing waves. We refer to comprehensive reviews, e.g., \cite{appelo:2009,astley:2000,gerdes:2000,rasch:1986} and references therein for a complete overview of the open boundary conditions issue.

In order to discretize differential problems on unbounded domains, spectral collocation methods \cite{shen:2001} based on scaled basis functions \cite{shen:2009b,shen:2009,wang:2009,xia:2021,zhuang:2010,zhuang:2010b} were used in unbounded regions coupled with finite volume \cite{benacchio:2013} and discontinuous Galerkin (DG) finite element \cite{benacchio:2019} methods on finite regions for one-dimensional linear and nonlinear hyperbolic problems. The coupled setup enabled accurate simulation of wavelike phenomena and economical damping of outgoing signals with minimal reflections in the region of interest using a small number of scaled Laguerre basis functions and tuning their scaling parameter in the semi-infinite part. The methodology was then broadened to include the diffusive case \cite{vismara:2022} within a seamless, extended DG (XDG) modal approach using Legendre basis functions in the bounded region, scaled spectral Laguerre basis functions in the unbounded region, and flux-based coupling of the two basis function sets at the finite/semi-infinite interface. A formal analysis confirmed the stability and the potential of the XDG framework in providing reliable and efficient numerical solutions in very large regions.

The present study extends the previous one-dimensional proof-of-concept models \cite{benacchio:2013,benacchio:2019,vismara:2022} in two unexplored directions:

\begin{itemize}
\item First, the XDG framework is implemented in two spatial dimensions and used to simulate linear and nonlinear problems on semi-infinite strip-shaped domains. The seamlessness of the approach is evident in the expressions for many of the integrals appearing in the discrete formulation in the finite and semi-infinite regions. Those expressions are formally identical up to replacing the finite-domain grid size with the inverse of the Laguerre scaling parameter;

\item Second, time-to-solution values with the two-dimensional XDG model are compared one-to-one with a single-domain DG framework that uses Legendre basis functions and an inhomogeneous grid with nodes coinciding with those of the extended scheme in the unbounded part. For several choices of model parameters, the new model provides numerical solutions with several times lower wallclock times at a given accuracy level, with the efficiency advantage growing with the number of Laguerre modes and thus larger simulated portions of the semi-infinite region.
\end{itemize}
We remark that alternative methods exist for high-order simulations on unbounded domains, but available methodologies are limited to the linear case without like-for-like efficiency comparisons \cite{kahana:2022}, concern elliptic problems \cite{gu:2021}, or lack formal stability analyses \cite{tatari:2014}. Recent studies\cite{chou:2022,xia:2021} introduced adaptive methods based on time-evolving distributions of Laguerre and Hermite collocation points for simulations on unbounded domains and applied them within spectrally adaptive neural network models \cite{xia:2022}. 

The rest of the paper is structured as follows.
Section \ref{sec: prob_def_geo_cont_wf} describes the problem definition and geometric setup, while section\ref{sec:discrete_formulation} provide details on the 2D XDG model and its discrete and algebraic formulation. Section \ref{sec:num_exp} reports the results obtained in (a) convergence tests with a standalone DG-Laguerre 2D formulation; (b) validation tests for the coupling strategy of the extended scheme for the linear advection-diffusion, nonlinear Burgers', and nonlinear shallow water equations; (c) performance tests in the semi-infinite regions; and (d) absorbing layer tests with an added reaction term modelling the damping of single-wave and wave train-like signals. Section \ref{sec:disc_conc} concludes the paper.

\section{Problem definition and geometry discretization}\label{sec: prob_def_geo_cont_wf}
In a two-dimensional semi-infinite strip $\Omega=[0,L_x]\times\mathbb{R}^+\equiv[0,+\infty)$, we consider the following conservation law with damping:
\begin{align}
    \begin{split}\label{eq:adv_eq_coupled}
    & \frac{\partial{q}}{\partial t}+\nabla\cdot\boldsymbol{\mathfrak{F}}(q,x,z)+\gamma(x,z)q=f \qquad (x,z)\in\Omega,\,t\in[0,T]\\
    & q(x,z,t)=q_{D}(x,z,t)\qquad\qquad\qquad\quad\textrm{if} \;(x,z)\in\partial\Omega_D \\
    & \boldsymbol{\mathfrak{F}}(q,x,z)\cdot\mathbf{n}=q_N(x,z,t)\qquad\qquad\quad\ \ \textrm{if} \;(x,z)\in\partial\Omega_N\\
    & \lim_{z\rightarrow+\infty}q(x,z,t)=0 \\
    & q(x,z,0) = q_0(x,z)
    \end{split}
\end{align}
where $\partial\Omega_D$ denotes the part of the boundary where a Dirichlet condition is imposed, while a Neumann condition is assigned to $\partial\Omega_N$. For the sake of simplicity, in the rest of this section we will focus on the transport-diffusion case, i.e.:

\begin{equation}
\boldsymbol{\mathfrak{F}}(q,x,z) = -\bm{\mu}\nabla q+\mathbf{F}   
\end{equation}

\noindent where $\bm{\mu}$ is a $2\times2$ symmetric and positive definite matrix uniformly bounded above and below, that is, there exist two positive constants $\mu_0$ and $\mu_1$ such that
\begin{equation}
    \mu_0\lVert w\rVert^2\leq w^\top\bm{\mu}w\leq \mu_1\lVert w\rVert^2 \qquad \forall w\in\mathbb{R}^2
\end{equation}
and we will assume that $\bm{\mu}$ is constant and diagonal:
\begin{equation}
    \bm{\mu}=\begin{bmatrix}\mu_x & 0 \\ 0 & \mu_z\end{bmatrix} \qquad \mu_x,\mu_z\in\mathbb{R}.
\end{equation}
If $\mu_x,\mu_z>0$, then $\bm{\mu}$ is symmetric and positive definite, and the boundedness condition is automatically satisfied with $\mu_0=\min(\mu_x,\mu_z)$ and $\mu_1=\max(\mu_x,\mu_z)$.
The function $\mathbf{F}$ denotes a generic transport term which may also depend on the spatial variables. Thus $\mathbf{F}=\mathbf{F}(q,x,z)=[F_1(q,x,z),F_2(q,x,z)]^\top$ but we will omit the dependence on $x$ and $z$ in the following expressions. We allow $\mathbf{F}$ to be non-linear in $q$ and discontinuous in $x$ and $z$. Hence, this formulation includes the linear advection-diffusion equation, for $\mathbf{F}(q,x,z)=\mathbf{u}(x,z)q$ with $\mathbf{u}=\mathbf{u}(x,z)=[u(x,z), v(x,z)]^\top$, and the Burgers' equation, for $\mathbf{F}(q,x,z)=[q^2/2,q^2/2]^\top$. The extension to the vectorial case and nonlinear systems is straightforward, see Section \ref{subsec: swe} below. Finally, the function $\gamma$ and the right-hand side $f=f(x,z,t)$ represent a reaction term and a given source/sink, respectively. 

The weak formulation of the problem is obtained by multiplying \eqref{eq:adv_eq_coupled} by a test function $\varphi\in V$ and integrating by parts:
For all test functions $\varphi\in V$, find $q\in V$ such that:
\begin{align}
\label{eq:adv_eq_gal_cont_coupled}
\begin{split}
 \frac{d}{dt}\int_\Omega q\varphi dx\,dz&-\int_{\partial\Omega}(\bm{\mu}\nabla q\cdot\mathbf{n}_{|_{\partial\Omega}})\varphi_{|_{\partial\Omega}} dx\,dz + \int_\Omega\bm{\mu}\nabla q\cdot\nabla\varphi dx\,dz+\\&+\int_{\partial\Omega}\mathbf{F}(q)\cdot\mathbf{n}{\varphi}_{|_{\partial\Omega}} dx\,dz-\int_\Omega \mathbf{F}(q)\cdot\nabla\varphi dx\,dz+\int_\Omega \gamma q\varphi dx\,dz=\int_\Omega f\varphi dx\,dz
 \end{split}
\end{align}
where $V$ is a sufficiently regular trial/test space.

We split the domain as $\Omega=([0,L_x]\times[0,L_z])\cup([0,L_x]\times[L_z,+\infty))$ and we discretize the problem using an XDG scheme\cite{vismara:2022}. 

In two dimensions, XDG is a coupled scheme consisting of:

\begin{itemize}
 \item a standard DG scheme in both directions in the rectangle $\Omega_{DG-DG}=[0,L_x]\times[0,L_z]$;
 \item a ``DG in $x$''-``Laguerre in $z$'' scheme with scaled Laguerre basis functions in the semi-infinite region $\Omega_{DG-LAG}=[0,L_x]\times[L_z,+\infty)$.
\end{itemize}

More specifically, we introduce in $\Omega_{DG-DG}$ the computational grid $\Omega_{DG-DG}=\bigcup_{m_x=1}^{N_x}\bigcup_{m_z=1}^{N_z}K_{m_x,m_z}$ consisting of $N_xN_z$ non-overlapping elements of area $\Delta x_{m_x}\Delta z_{m_z}$, where $K_{m_x,m_z}=[x_{m_x-1},x_{m_x}]\times[z_{m_z-1},z_{m_z}]$, so that $x_0=0$, $x_{N_x}=L_x$, $z_0=0$ and $z_{N_z}=L_z$. 

Denoting with $\mathbb{P}_{p_x}([x_{m_x-1},x_{m_x}])$ the space of polynomials of degree $p_x$ in $[x_{m_x-1},x_{m_x}]$, and with $\mathbb{P}_{p_z}([z_{m_z-1},z_{m_z}])$ the space of polynomials of degree $p_z$ in $[z_{m_z-1},z_{m_z}]$,  we define the space of polynomials of degree $p_x$ in $x$ and $p_z$ in $z$ on $K_{m_x,m_z}$ as:
\begin{equation}
 \mathbb{P}_{p_x,p_z}(K_{m_x,m_z})=\mathbb{P}_{p_x}([x_{m_x-1},x_{m_x}])\times\mathbb{P}_{p_z}([z_{m_z-1},z_{m_z}]).
\end{equation}
Then, we consider the discontinuous finite element space:
\begin{equation}
 V_h^{p_x,p_z}=\left\{v\in L^2(\Omega_{DG-DG}):\,v|_{K_{m_x,m_z}}\in\mathbb{P}_{p_x,p_z}(K_{m_x,m_z}),\;m_x=1,\ldots,N_x,\;m_z=1,\dots,N_z\right\}
\end{equation}
and choose as basis of $\mathbb{P}_{p_x,p_z}(K_{m_x,m_z})$ the normalized Legendre basis $\left\{\phi^{m_x}_j(x)\phi^{m_z}_i(z)\right\}_{j=0,\ldots,p_x,i=0,\ldots,p_z}$, with 
\begin{equation}
 \phi^{m_x}_j(x)=\sqrt{2j+1}L_j\left(2\frac{x-x_{m_x-1/2}}{\Delta x_{m_x}}\right) \qquad \phi^{m_z}_i(z)=\sqrt{2i+1}L_i\left(2\frac{z-z_{m_z-1/2}}{\Delta z_{m_z}}\right),
\end{equation}
where $L_j$ is the $j$-th Legendre polynomial. Hence, the solution of \eqref{eq:adv_eq_coupled} will be represented in each interval $K_{m_x,m_x}$ as
\begin{equation}\label{eq:solution_repr_DGDG}
    q(x,z,t)\approx\sum_{j=0}^{p_x}\sum_{i=0}^{p_z}q^{(j,i)}_{m_x,m_z}(t)\phi_j^{m_x}(x)\phi_i^{m_z}(z) \qquad (x,z)\in\Omega_{DG-DG}.
\end{equation}

Next, the semi-infinite region is discretized using the computational grid $\Omega_{DG-LAG}=\bigcup_{m_x=1}^{N_x}K_{m_x,\infty}$ where $K_{m_x,\infty}=[x_{m_x-1},x_{m_x}]\times[L_z,+\infty)$. Scaled Laguerre functions are chosen as basis functions in the $z$-direction \cite{benacchio:2013,benacchio:2019,vismara:2022,wang:2009}. Defining 
\begin{equation}
    \phi_i^\infty(z)=\hat{\mathscr{L}}_i^\beta(z-L_z) \qquad i=0,\dots,M,
\end{equation}
we can represent the solution as
\begin{equation}\label{eq:solution_repr_DGLAG}
    q(x,z,t)\approx\sum_{j=0}^{p_x}\sum_{i=0}^{M}q^{(j,i)}_{m_x,\infty}(t)\phi_j^{m_x}(x)\phi_i^{\infty}(z) \qquad (x,z)\in\Omega_{DG-LAG}.
\end{equation}
Note that the same notation $\phi$ with different superindices is used to denote both types of basis functions of the XDG approach in order to remark the unitary character of the discretization using Legendre-based and Laguerre-based elements.

\begin{figure}[ht]
\centering
\includegraphics{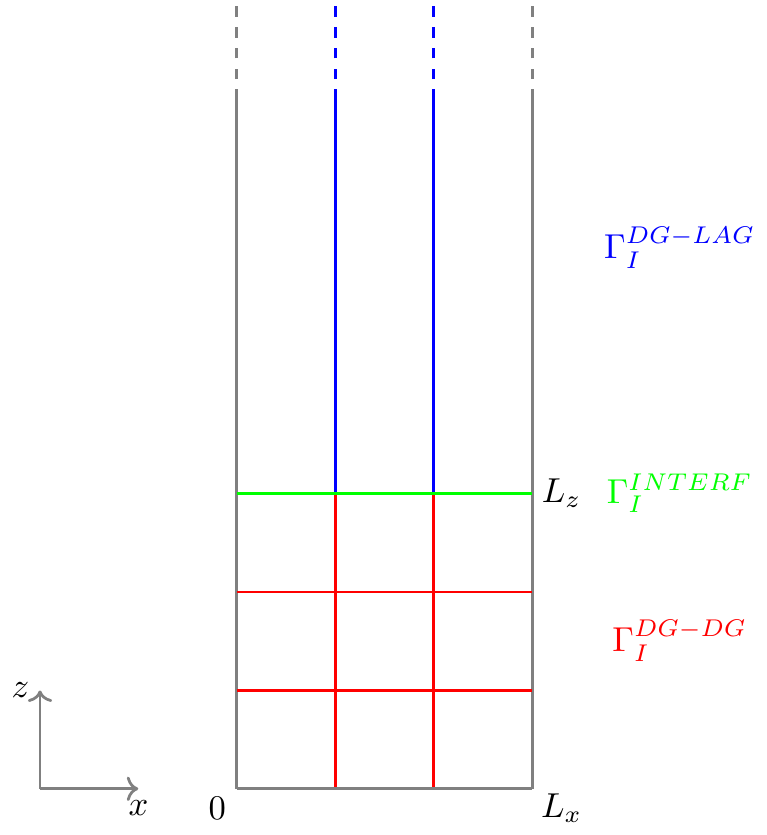}
\caption{The semi-infinite strip $\Omega$ with internal edges $\Gamma_I=\Gamma_I^{DG-DG}\cup\Gamma_I^{INTERF}\cup\Gamma_I^{DG-LAG}$, see text for details.}\label{fig:internal_edges_coupled}
\end{figure}

\noindent Written elementwise in $K_{m_x,m_z}$, $m_x=1,\dots,N_x$, $m_z=1,\dots,N_z,\infty$, \eqref{eq:adv_eq_gal_cont_coupled} reads
\begin{align}
\label{eq:adv_eq_gal_cont_loc_coupled}
\begin{split}
 &\frac{d}{dt}\int_{K_{m_x,m_z}} q\varphi dx\,dz-\int_{\partial{K_{m_x,m_z}}}(\bm{\mu}\nabla q\cdot\mathbf{n}_{|_{\partial{K_{m_x,m_z}}}})\varphi_{|_{\partial{K_{m_x,m_z}}}}\,ds + \int_{K_{m_x,m_z}}\bm{\mu}\nabla q\cdot\nabla\varphi dx\,dz+\\&+\int_{\partial{K_{m_x,m_z}}}\mathbf{F}(q)\cdot\mathbf{n}{\varphi}_{|_{\partial{K_{m_x,m_z}}}}\,ds-\int_{K_{m_x,m_z}} \mathbf{F}(q)\cdot\nabla\varphi dx\,dz+\int_{K_{m_x,m_z}}\gamma q\varphi dx\,dz=\int_{K_{m_x,m_z}} f\varphi dx\,dz.
 \end{split}
\end{align}
\noindent We then sum over all intervals $K_{m_x,m_z}$ to obtain
\begin{align}\label{eq:sum_over_intervals_cont_coupled}
\begin{split}
 \frac{d}{dt}\int_{\Omega} q\varphi dx\,dz&-\sum_{e\in\Gamma_h}\int_{e}(\bm{\mu}\nabla q\cdot\mathbf{n}_{|_e})\varphi_{|_e}\,ds + \int_{\Omega}\bm{\mu}\nabla q\cdot\nabla\varphi dx\,dz+\\&+\sum_{e\in\Gamma_h}\int_{e}\mathbf{F}(q)\cdot\mathbf{n}{\varphi}_{|_e}\,ds-\int_{\Omega} \mathbf{F}(q)\cdot\nabla\varphi dx\,dz+\int_{\Omega}\gamma q\varphi dx\,dz=\int_{\Omega} f\varphi dx\,dz
 \end{split}
\end{align}
where $\Gamma_h$ is the set of edges,
\begin{equation}
    \Gamma_h=\Gamma_I\cup\Gamma_D\cup\Gamma_N
\end{equation}
being $\Gamma_I$ the set of internal edges, and $\Gamma_D$ and $\Gamma_N$ the edges on $\partial\Omega$ where a Dirichlet or a Neumann boundary condition is imposed, respectively. We also allow periodic boundary conditions to be assigned on a subset of the boundary - in this case we will simply treat the corresponding edges as internal edges. The set $\Gamma_I$ is further subdivided as
\begin{equation}
    \Gamma_I = \Gamma_I^{DG-DG}\cup\Gamma_I^{DG-LAG}\cup\Gamma_I^{INTERF}
\end{equation}
to distinguish edges internal to $\Omega_{DG-DG}$, edges internal to $\Omega_{DG-LAG}$ and edges on the interface between the two regions (Figure \ref{fig:internal_edges_coupled}).

We analyze the second and fourth integral in \eqref{eq:sum_over_intervals_cont_coupled}, i.e., the integrals over the edges.
Starting from the diffusive term, we split the sum as
\begin{equation}
    -\sum_{e\in\Gamma_h}\int_{e}(\bm{\mu}\nabla q\cdot\mathbf{n}_{|_{e}})\varphi_{|_{e}} ds=-\sum_{e\in\Gamma_I\cup\Gamma_D}\int_{e}(\bm{\mu}\nabla q\cdot\mathbf{n}_{|_{e}})\varphi_{|_{e}} ds-\sum_{e\in\Gamma_N}\int_{e}(\bm{\mu}\nabla q\cdot\mathbf{n}_{|_{e}})\varphi_{|_{e}} ds.
\end{equation}
We impose the Neumann boundary condition on $\Gamma_N$, so that the second sum on the right hand side becomes
\begin{equation}
    -\sum_{e\in\Gamma_N}\int_eq_N\varphi_{|_e}ds.
\end{equation}
Regarding the first sum on the right hand side, we first notice that
\begin{equation}\label{eq:identity_jumps_avgs_coupled}
    -\sum_{e\in\Gamma_I\cup\Gamma_D}\int_{e}(\bm{\mu}\nabla q\cdot\mathbf{n}_{|_{e}})\varphi_{|_{e}} ds=-\sum_{e\in\Gamma_I\cup\Gamma_D}\int_e\{\{\bm{\mu}\nabla q\cdot\mathbf{n}\}\}[[\varphi]]ds -\sum_{e\in\Gamma_I\cup\Gamma_D}\int_e[[\bm{\mu}\nabla q\cdot\mathbf{n}]]\{\{\varphi\}\}ds,
\end{equation}
where the jump and average operators across an edge $e$ are defined as
\begin{equation}
    [[w]]=w^a-w^b \qquad \{\{w\}\}=\frac{w^a+w^b}{2},
\end{equation}
where the elements separated by $e$ are denoted as $\Omega^a$ and $\Omega^b$ and $\mathbf{n}$ is oriented from $\Omega^a$ to $\Omega^b$. If $e\in\Gamma_D$ is a boundary edge, we extend the definitions of jump and average as
\begin{equation}
    [[w]]=w^a \qquad \{\{w\}\}=w^a.
\end{equation}

\noindent If $q$ is regular, we have that $[[\bm{\mu}\nabla q\cdot\mathbf{n}]]=0$ and $[[q]]=0$ on every $e\in\Gamma_I\cup\Gamma_D$. Using this property, the Dirichlet boundary condition on $\Gamma_D$ on equation \eqref{eq:identity_jumps_avgs_coupled}, and the definition of jumps and averages on a boundary face, we obtain

\begin{equation}\label{eq:epsilon_sigma_coupled}
    \begin{split}
    -\sum_{e\in\Gamma_I\cup\Gamma_D}\int_{e}&(\bm{\mu}\nabla q\cdot\mathbf{n}_{|_{e}})\varphi_{|_{e}} ds=-\sum_{e\in\Gamma_I\cup\Gamma_D}\int_e\{\{\bm{\mu}\nabla q\cdot\mathbf{n}\}\}[[\varphi]]ds +\epsilon\sum_{e\in\Gamma_I\cup\Gamma_D}\int_e\{\{\bm{\mu}\nabla\varphi\cdot\mathbf{n}\}\}[[q]]ds+\\&\sum_{e\in\Gamma_I\cup\Gamma_D}\frac{\sigma}{\lvert e\rvert}\int_e[[q]][[\varphi]]ds + -\sum_{e\in\Gamma_D}\int_e\left(\epsilon\bm{\mu}\nabla\varphi_{|_e}\cdot\mathbf{n}+\frac{\sigma}{\lvert e\rvert}\varphi_{|_e}\right)q_D\,ds
    \end{split}
\end{equation}
for any choice of the penalization parameters $\epsilon$ and $\sigma$. $\lvert e\rvert$ denotes the length of the interface $e$.

We can now insert \eqref{eq:epsilon_sigma_coupled} into \eqref{eq:sum_over_intervals_cont_coupled} to obtain the following expanded weak formulation: for all test functions $\varphi\in V$, find $q\in V$ such that
\begin{multline}
\label{eq:weak_form_continuous_coupled}
 \frac{d}{dt}\int_{\Omega} q\varphi dx\,dz-\sum_{e\in\Gamma_I\cup\Gamma_D}\int_e\{\{\bm{\mu}\nabla q\cdot\mathbf{n}\}\}[[\varphi]]ds +\epsilon\sum_{e\in\Gamma_I\cup\Gamma_D}\int_e\{\{\bm{\mu}\nabla\varphi\cdot\mathbf{n}\}\}[[q]]ds+\\+\sum_{e\in\Gamma_I\cup\Gamma_D}\frac{\sigma}{\lvert e\rvert}\int_e[[q]][[\varphi]]ds + \int_{\Omega}\bm{\mu}\nabla q\cdot\nabla\varphi dx\,dz+\sum_{e\in\Gamma_h}\int_e\mathbf{F}(q)\cdot\mathbf{n}\varphi_{|_e}\,ds-\int_{\Omega} \mathbf{F}(q)\cdot\nabla\varphi dx\,dz+\int_\Omega\gamma q\varphi dx\,dz=\\=\int_{\Omega} f\varphi dx\,dz+\sum_{e\in\Gamma_D}\int_e\left(\epsilon\bm{\mu}\nabla\varphi_{|_e}\cdot\mathbf{n}+\frac{\sigma}{\lvert e\rvert}\varphi_{|_e}\right)q_D\,ds+\sum_{e\in\Gamma_N}\int_eq_N\varphi_{|_e}ds.
\end{multline}

\section{Algebraic formulation}\label{sec:discrete_formulation}

We can now move from the continuous to the discrete setting. We first introduce the spaces
\begin{align}
     &V_h^{p_x,p_z,M}=\left\{v\in V:\,v|_{\Omega_{DG-DG}}\in V_h^{p_x,p_z}, v|_{\Omega_{DG-LAG}}\in V_h^{p_x,M}\right\}\subset V, \\
     &V_h^{p_x,M}=\text{span}\{\phi_j^{m_x}\}_{j=0,\dots,p_x,m_x=1,\dots,N_x}\times\text{span}\{\phi_i^\infty\}_{i=0,\dots,M}.
\end{align}
When replacing $q\in V$ by $q_h\in V_h^{p_x,p_z,M}$ and $\varphi\in V$ by $\varphi_h\in V_h^{p_x,p_z,M}$ in the last expression in the previous section, the term $\mathbf{F}(q)\cdot\mathbf{n}{\varphi}_{|_e}$ is not well-defined, since functions in $V_h^{p_x,p_z,M}$ are not uniquely defined at the internal edges. We then replace it by $\hat{q}(q_h,\mathbf{n}){\varphi_h}_{|_e}$, where $\hat{q}$ is the numerical flux (see Appendix \ref{app: discrete_form} for details).

\noindent The Galerkin formulation associated with \eqref{eq:adv_eq_coupled}-\eqref{eq:weak_form_continuous_coupled} reads:

Find $q_h\in V_h^{p_x,p_z,M}$ such that, for all test functions $\varphi_h\in V_h^{p_x,p_z,M}$,
\begin{align}\label{eq:weak_form_final_disc_coupled}
    \begin{split}
 \frac{d}{dt}\int_{\Omega} q_h\varphi_h dx\,dz&-\sum_{e\in\Gamma_I\cup\Gamma_D}\int_e\{\{\bm{\mu}\nabla q_h\cdot\mathbf{n}\}\}[[\varphi_h]]ds +\epsilon\sum_{e\in\Gamma_I\cup\Gamma_D}\int_e\{\{\bm{\mu}\nabla \varphi_h\cdot\mathbf{n}\}\}[[q_h]]ds+\\&+\sum_{e\in\Gamma_I\cup\Gamma_D}\frac{\sigma}{\lvert e\rvert}\int_e[[q_h]][[\varphi_h]]ds + \int_{\Omega}\bm{\mu}\nabla q_h\cdot\nabla\varphi_h dx\,dz+\\&+\sum_{e\in\Gamma_h}\int_e\widehat{q}(q_h,\mathbf{n}){\varphi_h}_{|_e}\,ds-\int_{\Omega} \mathbf{F}(q_h)\cdot\nabla\varphi_h dx\,dz+\int_\Omega\gamma q_h\varphi_h dx\,dz=\int_{\Omega} f\varphi_h dx\,dz + \\&+\sum_{e\in\Gamma_D}\int_e\left(\epsilon\bm{\mu}\nabla{\varphi_h}_{|_e}\cdot\mathbf{n}+\frac{\sigma}{\lvert e\rvert}{\varphi_h}_{|_e}\right)q_D\,ds+\sum_{e\in\Gamma_N}\int_eq_N{\varphi_h}_{|_e}ds.
 \end{split}
\end{align}

\noindent In order to obtain a discrete formulation, the representations \eqref{eq:solution_repr_DGDG} or \eqref{eq:solution_repr_DGLAG} for $q_h$ are replaced in the above expression, choosing $\varphi_h=\phi_{j^\prime}^{m_x^\prime}\phi_{i^\prime}^{m_z^\prime}$ as test functions, with $j^\prime=0,\dots,p_x$, $m_x^\prime=1,\dots,N_x$, $m_z^\prime=1,\dots,N_z,\infty$ and $i^\prime=0,\dots,p_z$ if $m_z^\prime=1,\dots,N_z$, $i^\prime=0,\dots,M$ if $m_z^\prime=\infty$. The resulting integrals are then simplified using the orthogonality of the Legendre and Laguerre basis functions, see Appendix \ref{app: discrete_form} for details on the computation of the integrals.

In order to cast the spatially discrete problem in algebraic form, we first define the unknown vector as $\mathbf{q}\in\mathbb{R}^{N_{xtd}}$, with $N_{xtd}= N_xN_z(p_x+1)(p_z+1)+N_x(p_x+1)(M+1)$ . We also define $N_{bnd}=N_xN_z(p_x+1)(p_z+1)$ and $N_{unbnd}=N_x(p_x+1)(M+1)$ as the number of degrees of freedom in the bounded and unbounded region, respectively. We order the unknown degrees of freedom as follows. We first consider the degrees of freedom in $\Omega_{DG-DG}$ and define the local vector
\begin{equation}
    \mathbf{q}^{DG-DG}_{m_x,m_z}=[q_{m_x,m_z}^{(0,0)},\dots,q_{m_x,m_z}^{(p_x,0)},q_{m_x,m_z}^{(0,1)},\dots,q_{m_x,m_z}^{(p_x,1)},\dots,q_{m_x,m_z}^{(0,p_z)},\dots,q_{m_x,m_z}^{(p_x,p_z)}]^T\in\mathbb{R}^{(p_x+1)\times(p_z+1)},
\end{equation}
which contains the unknowns corresponding to $K_{m_x,m_z}$. The global vector for the DG-DG discretization in the finite region is then
\begin{equation}
    \mathbf{q}^{DG-DG}=[\mathbf{q}_{1,1},\dots,\mathbf{q}_{N_x,1},\mathbf{q}_{1,2},\dots,\mathbf{q}_{N_x,2},\dots,\mathbf{q}_{1,N_z},\dots,\mathbf{q}_{N_x,N_z}]^T\in\mathbb{R}^{N_{bnd}}.
\end{equation}
Concerning the DG-Laguerre discretization in the semi-infinite region, we order the unknowns as \begin{align}\label{eq:unknown_vector_DGLAG}
\begin{split}
    \mathbf{q}^{DG-LAG}=[q_{1}^{(0,0)},\dots,q_{1}^{(p_x,0)},q_{2}^{(0,0)},\dots,q_{2}^{(p_x,0)},\dots,q_{N_x}^{(0,0)},\dots,q_{N_x}^{(p_x,0)},\\
    q_{1}^{(0,1)},\dots,q_{1}^{(p_x,1)},q_{2}^{(0,1)},\dots,q_{2}^{(p_x,1)},\dots,q_{N_x}^{(0,1)},\dots,q_{N_x}^{(p_x,1)},\\\dots\\
    q_{1}^{(0,M)},\dots,q_{1}^{(p_x,M)},q_{2}^{(0,M)},\dots,q_{2}^{(p_x,M)},\dots,q_{N_x}^{(0,M)},\dots,q_{N_x}^{(p_x,M)}]^T\in\mathbb{R}^{N_{unbnd}}.
\end{split}
\end{align}
The global unknown vector for the XDG formulation will then be $\mathbf{q}=[\mathbf{q}^{DG-DG},\mathbf{q}^{DG-LAG}]^T\in\mathbb{R}^{N_{xtd}}$. The algebraic formulation consists of the following components:

\begin{itemize}
\item The equation for the time derivative of the two unknown vectors is (see also equation \eqref{eq:1st_coupled}):
\begin{equation}
    \Delta x\Delta z\frac{d\mathbf{q}^{DG-DG}}{dt} \qquad \frac{\Delta x}{\beta}\frac{d\mathbf{q}^{DG-LAG}}{dt}.
\end{equation}

\item The terms referring to the degrees of freedom in the bounded region result in the algebraic representation $A^{DG-DG}\mathbf{q}^{DG-DG}$, where $A^{DG-DG}\in\mathbb{R}^{N_{bnd}\times N_{bnd}}$ represents the DG-DG discretization in $\Omega_{DG-DG}$ (see expressions \eqref{eq:diff_internal_DGDG_vertical_coupled},\eqref{eq:diff_internal_DGDG_horizontal_coupled},\eqref{eq:diff_dirichlet_coupled_left_DGDG},\eqref{eq:diff_dirichlet_coupled_right_DGDG},\eqref{eq:diff_dirichlet_coupled_bottom} and the first terms of equation \eqref{eq:diff_interface_DGLAG_coupled}). 

\item The terms referring to the degrees of freedom in the unbounded region result in the algebraic formulation $A^{DG-LAG}\mathbf{q}^{DG-LAG}$, where $A^{DG-LAG}\in\mathbb{R}^{N_{unbnd}\times N_{unbnd}}$ stems from the DG-Laguerre discretization in $\Omega_{DG-LAG}$ (see expressions \eqref{eq:diff_internal_DGLAG_coupled},\eqref{eq:diff_dirichlet_coupled_left_DGDG},\eqref{eq:diff_dirichlet_coupled_right_DGDG} with $m_z=\infty$ and $\Delta z =1/\beta$ and the fourth term of equation \eqref{eq:diff_interface_DGLAG_coupled}). 

\item The coupling terms are expressed algebraically  as $A^{C,12}\mathbf{q}^{DG-DG}$ and $A^{C,21}\mathbf{q}^{DG-LAG}$ (see the second term and the third term of equation \eqref{eq:diff_interface_DGLAG_coupled}). The matrices
\begin{equation}
    A^{C,12}\in\mathbb{R}^{N_{unbnd}\times N_{bnd}}\quad \text{ and }\quad A^{C,21}\in\mathbb{R}^{N_{bnd}\times N_{unbnd}}
\end{equation}
are the coupling matrices, which are non-zero only in the last $N_x(p_x+1)(p_z+1)$ columns and the last $N_x(p_x+1)(p_z+1)$ rows, respectively, corresponding to the degrees of freedom in the $N_x$ DG-DG elements adjacent to the interface. 

\item The transport term can be assembled in the vector $\mathbf{F}(\mathbf{q})\in\mathbb{R}^{N_{xtd}}$ (see expressions \eqref{eq:adv_ver_interfaces} and \eqref{eq:adv_hor_interfaces} and Appendix \ref{sec: appendix_flux_int}). 

\item The damping term results in the algebraic terms $G^{DG-DG}\mathbf{q}$ and $G^{DG-LAG}\mathbf{q}$, where $G^{DG-DG}\in\mathbb{R}^{N_{bnd}\times N_{bnd}}$ and $G^{DG-LAG}\mathbb{R}^{N_{unbnd}\times N_{unbnd}}$ (from \eqref{eq:reac_DGDG_coupled} and \eqref{eq:reac_DGLAG_coupled}). 

\item Finally, the right-hand side vector is computed as the vector $\mathbf{f}\in\mathbb{R}^{N_{xtd}}$ (from equations \eqref{eq:rhs_DGDG_coupled_1st},\eqref{eq:rhs_DGDG_coupled_2nd},\eqref{eq:rhs_DGDG_coupled_3rd},and \eqref{eq:rhs_DGLAG_coupled_2nd},\eqref{eq:rhs_DGLAG_coupled_3rd}).

\end{itemize}

\noindent The algebraic formulation of the problem is then
\begin{equation}\label{eq:adv_diff_reac_algebraic}
    M\frac{d\mathbf{q}(t)}{dt}+A\mathbf{q}(t)+\mathbf{F}(\mathbf{q})+G\mathbf{q}(t)=\mathbf{f}(t),
\end{equation}
where
\begin{equation}
    M=\begin{bmatrix}\Delta x\Delta zI & 0 \\ 0 & \frac{\Delta x}{\beta}I\end{bmatrix} \qquad A=\begin{bmatrix}A^{DG-DG} & A^{C,21} \\ A^{C,12} & A^{DG-LAG}\end{bmatrix} \qquad G=\begin{bmatrix}G^{DG-DG} & 0 \\ 0 & G^{DG-LAG}\end{bmatrix}.
\end{equation}
In the linear case, the transport term reduces to $F\mathbf{q}$, where the matrix $F$ has the same block structure as $A$.

The spatially discretized problem can be solved in time by means of a suitable time integration scheme. In the linear case, equation \eqref{eq:adv_diff_reac_algebraic} reads
\begin{equation}
    M\frac{d\mathbf{q}}{dt}+B\mathbf{q}=\mathbf{f},
\end{equation}
where $B=A+F+G$. The $\theta$-method can then be employed: given $\mathbf{q}^0$ and $\theta\in[0,1]$, the solution at time step $n$ is given by
\begin{equation}
    (M+\Delta t\theta B)\mathbf{q}^{n+1}=[M-\Delta t(1-\theta)B]\mathbf{q}^n+\Delta t\theta\mathbf{f}^{n+1}+\Delta t(1-\theta)\mathbf{f}^n.
\end{equation}

\noindent In the non-linear case we rewrite \eqref{eq:adv_diff_reac_algebraic} as
\begin{equation}\label{eq:semidisc_form_for_imex}
    \frac{d\mathbf{q}}{dt}=-M^{-1}(A+G)\mathbf{q}-M^{-1}(\mathbf{F}(\mathbf{q})-\mathbf{f}),
\end{equation}
then separate the linear and non-linear terms by defining $\mathbf{f}_{l}(\mathbf{q},t)=-M^{-1}(A+G)\mathbf{q}$ and $\mathbf{f}_{nl}(\mathbf{q},t)=-M^{-1}(\mathbf{F}(\mathbf{q})-\mathbf{f})$. Problem \eqref{eq:semidisc_form_for_imex} is discretized in time using a $3$-stage IMEX-ARK method \cite{giraldo:2013}.

The algebraic formulation readily extends to the vectorial case. The theoretical details are omitted here for conciseness, section \ref{subsec: swe} contains numerical results with the shallow water equations.

\section{Numerical experiments}\label{sec:num_exp}
In this section, several numerical tests of the XDG scheme are carried out. First, we perform a convergence analysis of the standalone unbounded two-dimensional section of the XDG scheme, where the $x$ direction is discretized with Legendre basis functions and the $z$ direction is discretized with scaled Laguerre basis functions. Therefore, we set $L_z=0$ and compute the convergence rates of both the spatial discretization ($x$-convergence and $z$-convergence) and the temporal discretization. Next, the coupling strategy developed in Section \ref{sec:discrete_formulation} between the finite $\Omega_{DG-DG}$ and the semi-infinite $\Omega_{DG-LAG}$ portions of the 2D strip in the XDG model is validated both on the linear 2D advection-diffusion equation and on non-linear problems - the 2D Burgers' equation and the 2D shallow water equations. As in previous work \cite{benacchio:2013,benacchio:2019,vismara:2022}, the quality of the coupling is evaluated by computing the errors of the XDG approach with respect to a reference DG discretization on a larger domain in benchmarks with Gaussian signals. Next, we compare the performance of both the 2D XDG scheme and a DG discretization that uses a non-uniform grid in the semi-infinite direction. For both schemes, we compute wallclock times and errors with respect to an exact solution in tests with the advection-diffusion equation and Gaussian initial data using a large number of Laguerre modes and small values of the scaling parameter so as to simulate a larger portion of the unbounded region. Finally, an absorbing layer is implemented in the semi-infinite region by means of a sigmoidal reaction term, and tests are run with wave trains to check whether the obtained model is able to efficiently damp outgoing signals without reflections in the finite part of the two-dimensional strip.  

\noindent Errors are computed using suitable Gaussian quadrature rules. Specifically, we define the discrete norms
\begin{align}
    &\lVert q_h\rVert_{L^2}=\sqrt{\sum_{m_x=1}^{N_x}\sum_{m_z=1}^{N_z}\displaystyle\frac{\Delta x_{m_x}}{2}\displaystyle\frac{\Delta z_{m_z}}{2}\sum_{k=1}^{ng_x}\sum_{l=1}^{ng_z}\left[q_h\left(x_{m_x}+\displaystyle\frac{\Delta x_{m_x}}{2}(\hat{x}_k-1),z_{m_z}+\displaystyle\frac{\Delta z_{m_z}}{2}(\hat{z}_l-1)\right)\right]^2w^x_kw^z_l} \\
    &\lVert q_h\rVert_{L^\infty}=\max_{m_x=1,\dots,N_x}\max_{m_z=1,\dots,N_z}\max_{k=1,\dots,ng_x}\max_{l=1,\dots,ng_z}\left\lvert q_h\left(x_{m_x}+\displaystyle\frac{\Delta x_{m_x}}{2}(\hat{x}_k-1),z_{m_z}+\displaystyle\frac{\Delta z_{m_z}}{2}(\hat{z}_l-1)\right)\right\rvert
\end{align}
where $\{\hat{x}_k,w_k^x\}_{k=1}^{ng_x}$ and $\{\hat{z}_l,w_l^z\}_{l=1}^{ng_z}$ are Gaussian nodes and weights in the reference interval $[-1,1]$. Then, we define the absolute and relative errors between the numerical solution $q_h$ and a reference solution $q_{ref}$ as
\begin{equation}
    \mathcal{E}_s^{abs}:=\lVert q_h-q_{ref}\rVert_{L^s} \qquad \mathcal{E}_s^{rel}:=\displaystyle\frac{\lVert q_h-q_{ref}\rVert_{L^s}}{\lVert q_{ref}\rVert_{L^s}} \qquad s=2,\infty.
\end{equation}
In all the following tests we employ a uniform grid in time with $N_t$ time steps, so that $\Delta t=T/N_t$. We compute absolute errors at the final time and we estimate the rates of convergence as follows. Assume that $\mathcal{E}_s(h)\approx Ch^{r_s}$, where $C$ is a constant and $h=\Delta x$ or $h=\Delta t$. Then
\begin{equation}
    \displaystyle\frac{\mathcal{E}_s(h_1)}{\mathcal{E}_s(h_2)}\approx \left(\frac{h_1}{h_2}\right)^{r_s}
\end{equation}
so that
\begin{equation}
    {r_s}\approx\log_{\frac{h_1}{h_2}}\left(\displaystyle\frac{\mathcal{E}_s(h_1)}{\mathcal{E}_s(h_2)}\right) \qquad s=2,\infty.
\end{equation}
The Courant number values in the $x$ and $z$ direction are computed as $C_x=u_x\Delta tp_x/\Delta x$ and $C_z=u_z\Delta tp_z/\Delta z$; in the semi-infinite region, $\Delta z$ is taken as the distance between the first two Laguerre nodes and $p_z=1$. 
In all experiments on the linear advection-diffusion equation we employ the Crank-Nicolson method in time, that is, the $\theta$-method with $\theta=1/2$, while the Burgers' equation and the non-linear shallow water equations are solved in time by means of a $3$-stage IMEX-ARK method \cite{giraldo:2013}. Concerning the DG penalization method, we focus on NIPG ($\epsilon=1$), which is stable for any choice of $\sigma\geq0$; we therefore set $\sigma=0$ in the following.
For the sake of conciseness of the main text, detailed results are placed in the Appendix, see the tables in Appendix \ref{sec:tables}. The model was implemented in MATLAB$^\textrm{\textregistered}$ Version 9.9.0.1524771 (R2020b Update 2) on a workstation with processor CPU Intel$^\textrm{\textregistered}$ Core$\textsuperscript{TM}$ i7-9750H CPU @ 2.60GHz with 6 cores and 32GB RAM. Timings reported in the performance tests in section \ref{subsec: undamped_runs_si_domain} are therefore to be taken as preliminary. An implementation with a lower-level programming language will allow more accurate performance measurements via batch submissions to compute-only nodes on HPC systems - this is left for future work. 

\subsection{Standalone DG-Laguerre discretization - Convergence tests}\label{subsec:conv_tests}
The first set of numerical experiments of the two-dimensional model focuses on a single-domain discretization of the semi-infinite strip, with a Laguerre discretization in the $z$ direction and a DG scheme in the $x$ direction ($L_z=0$). A convergence analysis is carried out on the linear undamped non-homogeneous 2D advection-diffusion equation, i.e., equation \eqref{eq:adv_eq_coupled} with $\mathbf{F}(q,x,z)=[u_x, u_z]^\top q(x,z,t)$, $u_x,u_z\in\mathbb{R}$, $\gamma=0$. We impose the exact solution:
\begin{equation}
    q(x,z,t)=\text{exp}\left[-\left(\displaystyle\frac{x-x_0}{\sigma_0}\right)^2\right]ze^{-z}\text{sin}^2(z-t)
\end{equation}
with $x_0=L_x/2$ and $\sigma_0=L_x/10$, and compute the right-hand side accordingly. We impose homogeneous Dirichlet boundary conditions on the lower boundary $z=0$ and periodic conditions on the right and left boundaries $x=0$ and $x=L_x$. Absolute errors are computed at the final time $T$ with respect to the exact solution.

First, we evaluate spatial convergence in the $z$ direction. To this end, a large number of intervals are used for the DG discretization in the $x$ direction, $N_x=1000$ with $p_x=3$. In order to highlight the convergence in $z$, we also set a small time step, $\Delta t=5\times 10^{-4}$, and we run the simulation for $N_t=100$ time steps. We choose $M\in\{5,10,20,30,35\}$; the corresponding values of the Courant number in the $z$-direction are $C_z\in\{8.10\times10^{-3},1.40\times10^{-2},2.86\times10^{-2},4.22\times10^{-2},4.90\times10^{-2}\}$, while $C_x=1.5$. The small values of the Courant number in the $z$ direction are due to the time step $\Delta t$ being much smaller than the grid spacing $\Delta z$. This is not dictated by stability reasons, as the implicit time integration scheme is stable for values of $C_z=O(1)$ (this is confirmed by numerical tests not reported here). Our intent is rather to highlight the convergence in $z$ by making the discretization error in time negligible. Experimental convergence is exponential in the number of Laguerre modes $M$ (straight lines in the semi-logarithmic plot, Figure \ref{fig:conv_tests} left, see also table \ref{tab:M_N300} in Appendix \ref{sec:tables}), in line with expectations with a spectral discretization. In Figure \ref{fig:conv_tests} center, the coefficient of the exponential (solid black line) was chosen as $0.4$ for benchmarking purposes. Behaviour for larger $M$ values (not shown) confirms the super-polynomial convergence. 

Next, we study the convergence of the DG discretization in the $x$ direction. The NIPG scheme is stable for any choice of $\sigma\geq0$, but it is suboptimal when the polynomial degree is even, in agreement with the theory \cite{riviere:2008}. Therefore, we get quadratic convergence for $p_x=2$ as well as $p_x=1$, but convergence of order $4$ is recovered when $p_x=3$ (Figure \ref{fig:conv_tests} center for $p_x=3$, see also Table \ref{tab:x_err_eps1} in appendix \ref{sec:tables} for the results for all $p_x=1,2,3$). The time step and the number of time steps are the same as the previous test, ($\Delta t=5\times10^{-4}\,s$, $N_t=100$).

Finally, for a space-time convergence test, we choose $M=60$ and fix the horizontal Courant number $C_x=0.5$, with $p_x=1$. We then refine $\Delta t$ and $\Delta x$ simultaneously at constant $C_x$, considering the pairs \\ $(N_x,N_t) \in \{(50,100), (100,200), (150,300), (200,400), (250,500)\}$. The corresponding values of $C_z$ are $C_z\in\{1.99,1.00,0.66,0.50,0.40\}$.  Errors converge quadratically with increasing resolution (Figure \ref{fig:conv_tests} right, see also Table \ref{tab:x_t_conv_eps1} in Appendix \ref{sec:tables}).

\begin{figure}[tb!]
\centering
\includegraphics{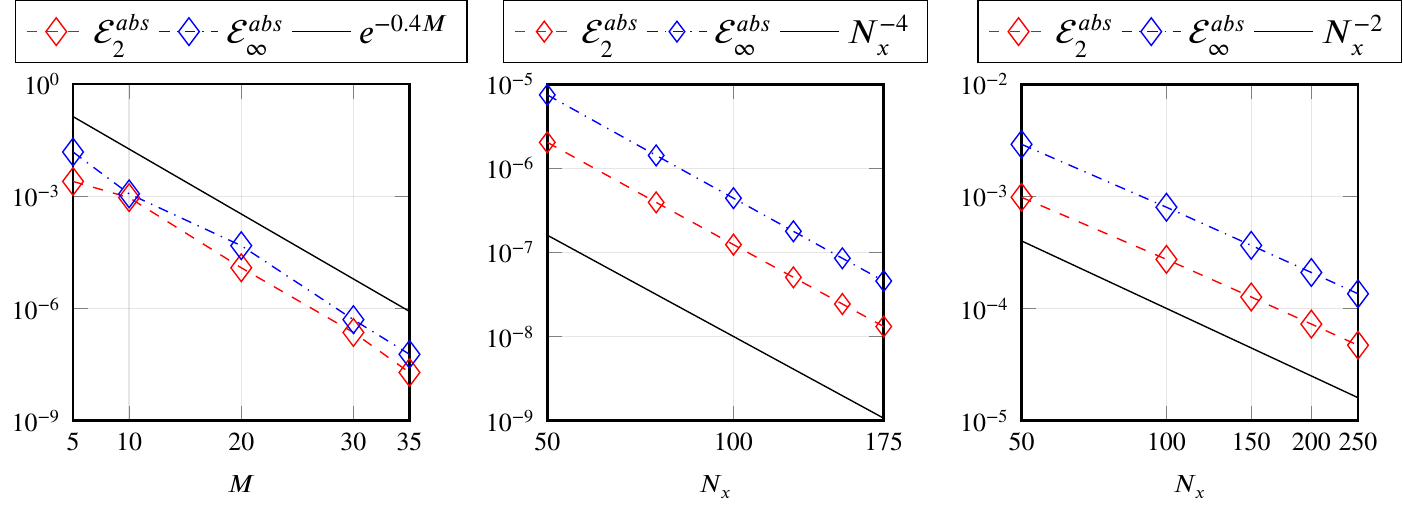}
\caption{Convergence tests of the 2D XDG scheme, linear advection-diffusion equation on a semi-infinite strip. Space convergence in $M$ (left),  space convergence in $N_x$, $p_x=3$ (center), space-time convergence, $p_x=1$ (right). For all plots, $L_x=1\,m$,$\mu_x=0.05\,m^2/s$, $\mu_z=0.01\,m^2/s$, $u_x=1\,m/s$, $u_z=2\,m/s$. For the left plot,  $N_x=1000$, $p_x=3$, $\beta=5$, $T=5\times 10^{-2}\,s$, $N_t=100$. For the center plot, $p_x=3$, $\beta=6$, $M=60$,  $T=5\times10^{-2}\,s$, $N_t=100$. For the right plot, $p_x=1$, $\beta=6$, $M=60$, $T=1\,s$, $N_t=2N_x$.}\label{fig:conv_tests}
\end{figure}

\subsection{Validation of the coupling strategy in the 2D XDG scheme}\label{subsec:coupling_valid}
In the second set of tests, we validate the accuracy of the coupling between the finite part and the semi-infinite part of the 2D strip, comparing the results obtained with the 2D XDG scheme and those obtained with a standalone DG discretization on a larger domain \cite{benacchio:2013,benacchio:2019,vismara:2022}, still considering the undamped case $\gamma=0$. We remark that the current implementation of the XDG model enables the use of polynomial degree up to $4$ for the finite part, as shown in the convergence tests and later in the efficiency tests of section \ref{subsec: undamped_runs_si_domain}. However, since a high polynomial degree is not needed to validate the coupling strategy, we restrict ourselves to $p_x=p_z=1$ in the coupling validation tests.

\subsubsection{Linear advection-diffusion equation}
\begin{figure}[ht]
\centering
    \includegraphics[width=0.95\textwidth]{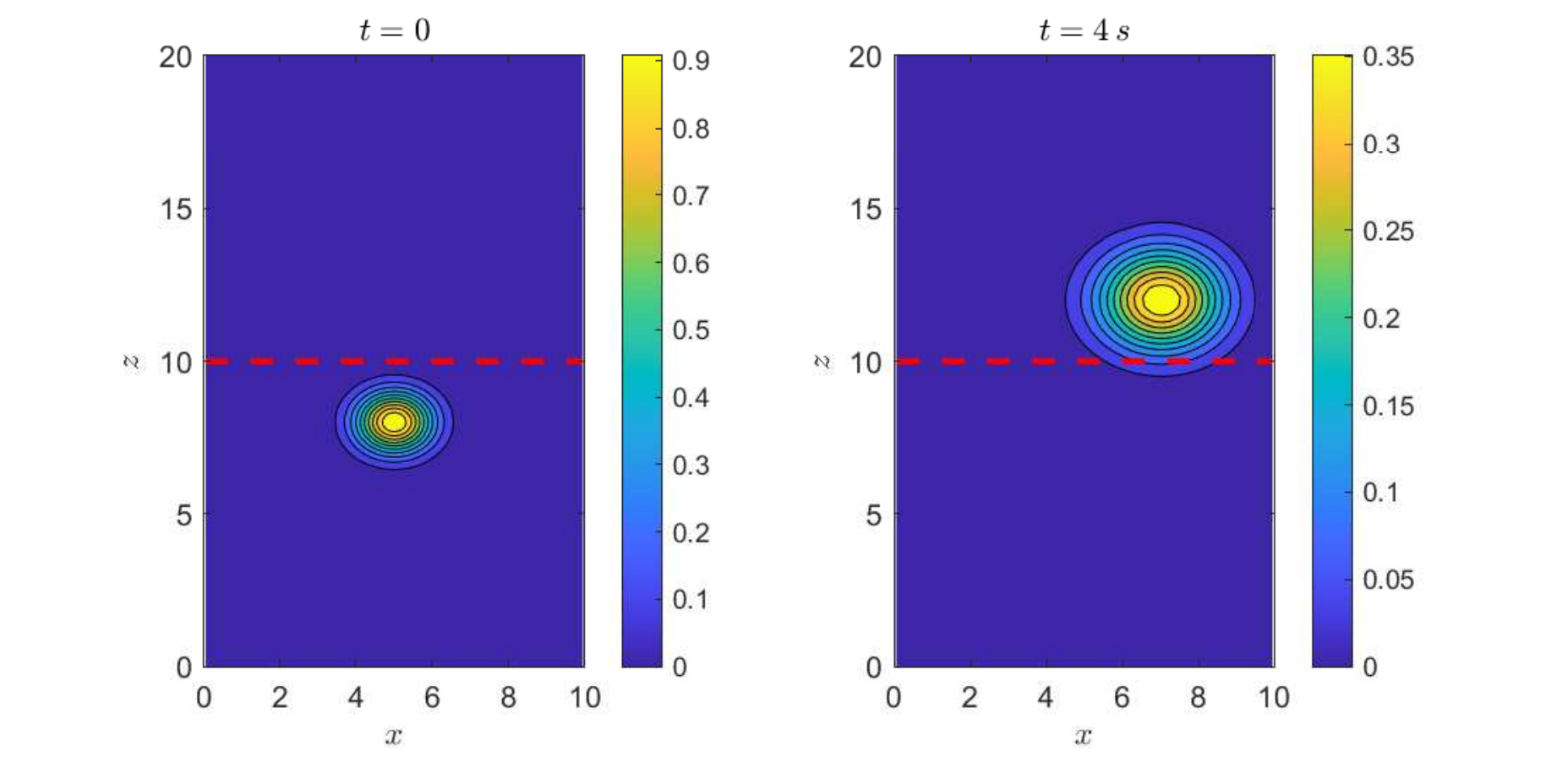}
\caption{ 2D Advection-diffusion equation. Initial profile (left, contour interval 0.1) and solution with the XDG scheme at the final time $T=4\,s$ (right, contour interval 0.05). The red dashed line denotes the interface between finite and semi-infinite regions. $L_x=L_z=10\,m$, $p_x=p_z=1$, $N_x=50$, $N_z=500$, $M=40$, $\beta=6$, $N_t=200$, $\mu_x=\mu_z=0.1\,m^2/s$, $u_x=0.5\,m/s$, $u_z=1\,m/s$, $\sigma_x=\sigma_z=1\,m$.}\label{fig:validation_linear_ux05}
\end{figure}

 We start with the advection-diffusion equation with Gaussian initial data placed inside the finite part.  More specifically, we set $L_x=L_z=10\,m$, $p_x=p_z=1$, $N_x=50$, $N_z=500$, $T=4\,s$, $N_t=200$, corresponding to $\Delta x=0.2\,m$, $\Delta z=0.02\,m$, $\Delta t=0.02\,s$. The physical parameters are $\mu_x=\mu_z=0.1\,m^2/s$, $u_x=0.5\,m/s$, $u_z=1\,m/s$, so that the Courant number in the $z$-direction is $C_z=1$. A homogeneous Dirichlet condition is imposed on the lower boundary, with periodic conditions on the vertical boundaries. The initial data is the Gaussian profile
\begin{equation}\label{eq:gauss_initcond}
    q_0(x,z)=A\text{exp}\left[-\left(\frac{x-x_0}{\sigma_x}\right)^2\right]\text{exp}\left[-\left(\frac{z-z_0}{\sigma_z}\right)^2\right]
\end{equation}
with $A=1\,m$, $x_0=5\,m$ and $z_0=8\,m$. As the simulation progresses, the initial profile is transported and diffused and crosses the finite/semi-infinite interface. We test the coupling approach for different choices of $\sigma_x=\sigma_z$ and $M$. For a fixed value of $M$, the scaling parameter $\beta$ should be chosen carefully in order to represent the semi-infinite part of the domain taking care of two aspects. On the one hand, the spacing between the first two Laguerre nodes should be comparable to the size of the last DG interval for the sake of information exchange at the interface. On the other hand, no perturbation leaving the finite domain should reach the last Laguerre node, so $\beta$ should be tuned in such a way that the semi-infinite part of the domain is large enough for the problem at hand. 
The solution computed by the 2D XDG model  (Figure \ref{fig:validation_linear_ux05}) with the combinations $(M,\beta)\in\{(10,4),(40,6)\}$ is compared to a reference 2D single-domain DG discretization on $[0,L_x]\times[0,2L_z]$ with the same spacing in the $z$-direction as the finite region of the extended scheme, that is, $N_z^\prime=1000$. Relative errors of the XDG scheme with respect to the reference are at most around a few percent (Table \ref{tab:adv_diff_coup_valid} in Appendix \ref{sec:tables}), in line with previous published work \cite{vismara:2022}. It is also worth noting that, with the current set of parameters, the number of total entries of the matrix of the XDG scheme is almost $4$ times smaller than the corresponding matrix for the standalone DG scheme. For $M=10$ and $M=40$ Laguerre basis functions, the matrix of the XDG scheme has about half as many non-zero entries as the matrix for the standalone DG scheme.

\subsubsection{Burgers' equation}
\begin{figure}[ht]
    \centering
    \includegraphics[width=\textwidth]{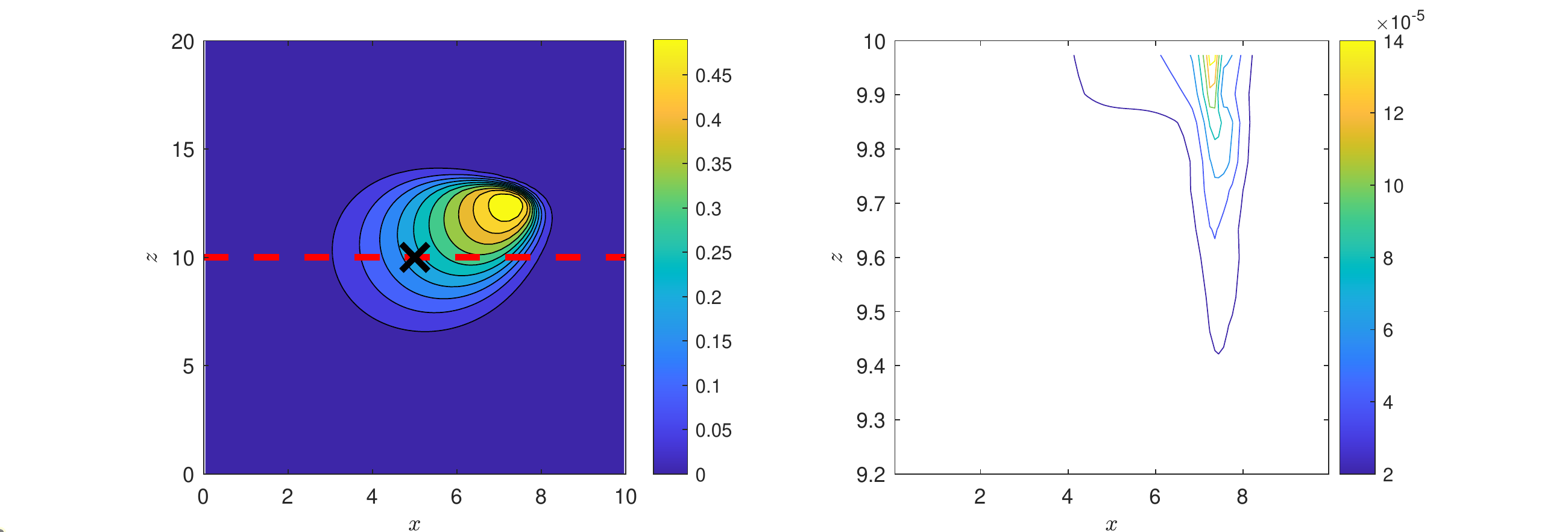}
    \caption{ Burgers' equation, Gaussian initial data. Contour plot of the numerical solution at the final time (left, contour interval 0.05) and absolute error with respect to a single-domain DG discretization (right, contour interval $2\times10^{-5}$, note the different range of values on the $z$ axis). The red dashed line denotes the interface between finite and semi-infinite regions. The black cross marker denotes the center of the initial profile. $L_x=L_z=10\,m$, $p_x=p_z=1$, $N_x=50$, $N_z=80$, $M=80$, $\beta=10$, $T=5\,s$, $N_t=500$, $\mu_x=\mu_z=0.05\,m^2/s$.}
    \label{fig:plot_burgers}
\end{figure}

Next, we consider the non-linear case and the Burgers' equation by setting $\mathbf{F}(q,x,z)=[q^2/2,q^2/2]^\top$ in \eqref{eq:adv_eq_coupled}, and the parameters $L_x=L_z=10\,m$, $N_x=50$, $N_z=80$, $\Delta x=0.2\,m$, $\Delta z=0.125\,m$, $p_x=p_z=1$, $\mu_x=\mu_z=0.05\,m^2/s$. We also set $\beta=10$. The initial condition is a Gaussian profile centered at the interface, that is, \eqref{eq:gauss_initcond} with $x_0=L_x/2=5\,m$ and $z_0=L_z=10\,m$. We also choose $\sigma_x=1\,m$ and $\sigma_z=2\,m$. The system is solved until the final time $T=5\,s$ with $N_t=500$ time steps. As above, we compare the solution of the XDG scheme (Figure \ref{fig:plot_burgers} left) with a reference standalone DG discretization on $[0,L_x]\times[0,2L_z]$ with the same grid spacing as in the bounded region of the XDG scheme, i.e., with $N_x^\prime=N_x=50$ and $N_z^\prime=2N_z=160$. Relative and absolute errors in the finite domain due to the coupling are small and limited to a narrow region near the interface (Figure \ref{fig:plot_burgers} right, see also Table \ref{tab:adv_diff_coup_valid} in Appendix \ref{sec:tables}).

\subsubsection{Non-linear shallow water equations}\label{subsec: swe}
The two-dimensional non-linear shallow water equations are \cite{leveque:2002}
\begin{equation}
    \begin{cases}
        \dpart{h}{t}+\dpart{(uh)}{x}+\dpart{(vh)}{z}=0 \\
        \dpart{(uh)}{t}+\dpart{(u^2h+\displaystyle\frac{gh^2}{2})}{x}+\dpart{(uvh)}{z}=0 \\
        \dpart{(vh)}{t}+\dpart{(uvh)}{x}+\dpart{(v^2h+\displaystyle\frac{gh^2}{2})}{z}=0
    \end{cases}
    \label{eq:nl_swe}
\end{equation}
where $h=h(x,z,t)$ is the water depth, $u=u(x,z,t)$ and $v=v(x,z,t)$ the velocities in the $x$ and $z$ direction, respectively, and $g$ is the acceleration of gravity. Introducing the vectors
\begin{equation}
    \mathbf{q}=\mathbf{q}(x,z,t)=\begin{bmatrix}
        q_1(x,z,t) \\ q_2(x,z,t) \\ q_3(x,z,t)
    \end{bmatrix}=\begin{bmatrix}
        h(x,z,t) \\ h(x,z,t)u(x,z,t) \\ h(x,z,t)v(x,z,t)
    \end{bmatrix}\in\mathbb{R}^3
\end{equation}
and
\begin{equation}
    \mathbf{E}(\mathbf{q})=\begin{bmatrix}
        uh \\ u^2h+gh^2/2 \\ uvh
    \end{bmatrix}=\begin{bmatrix}
        q_2 \\ q_2^2/q_1+gq_1^2/2 \\ q_2q_3/q_1
    \end{bmatrix} \qquad \mathbf{G}(\mathbf{q})=\begin{bmatrix}
        vh \\ uvh \\ v^2h+gh^2/2
    \end{bmatrix}=\begin{bmatrix}
        q_3 \\ q_2q_3/q_1 \\ q_3^2/q_1+gq_1^2/2
    \end{bmatrix}
\end{equation}
we obtain the system
\begin{equation}
    \displaystyle\frac{\partial\mathbf{q}}{\partial t}+\displaystyle\frac{\partial\mathbf{E}(\mathbf{q})}{\partial x}+\displaystyle\frac{\partial\mathbf{G}(\mathbf{q})}{\partial z}=0.
\end{equation}
We solve the problem in the semi-infinite strip $\Omega=[0,10\,m]\times[0,+\infty)$. To that end, we set $N_x=50$ and $p_x=1$; we place the interface at $L_z=8\,m$ with $N_z=40$ and $p_z=1$. We then test the accuracy of the XDG scheme for different choices of $M$ and $\beta$: the value of $\beta$ is selected so that the extension of the semi-infinite part is approximately the same in all tests. The final time is $T=0.75\,s$ with $N_t=750$ time steps. We write depth and velocities as $h(x,z,t)=H+\Tilde{h}(x,z,t)$, $u(x,z,t)=U+\Tilde{u}(x,z,t)$ and $v(x,z,t)=V+\Tilde{v}(x,z,t)$, where $\Tilde{h}$, $\Tilde{u}$ and $\Tilde{v}$ are perturbations of the reference configuration, taken to be $H=10\,m$, $U=1\,m/s$, $V=2\,m/s$. The initial condition is homogeneous for both $\Tilde{u}$ and $\Tilde{v}$. We consider for $\Tilde{h}(x,z,0)$ a Gaussian profile as in \eqref{eq:gauss_initcond} (see Figure \ref{fig:validation_linear_ux05} left) with $A=1\,m$, $x_0=L_x/2=5\,m$, $\sigma_x=\sigma_z=1\,m$; in the outgoing case, the initial profile is placed inside the finite region, with $z_0=5\,m$, while in the incoming case it is centered beyond the interface, at $z_0=10\,m$. The solution of the XDG scheme is compared to a reference single-domain DG solution on $L_z^\prime=3L_z=24\,m$ with the same grid spacing, $N_z^\prime=3N_z=120$\cite{benacchio:2013}. The boundary conditions are periodic in the $x$ direction and an outflow condition is imposed on the lower boundary $z=0$. With these setup choices, in both the incoming case and the outgoing case, the initial data is horizontally propagated across the periodic boundaries. 

\begin{figure}[htb!]
    \centering
    \includegraphics[width=\textwidth]{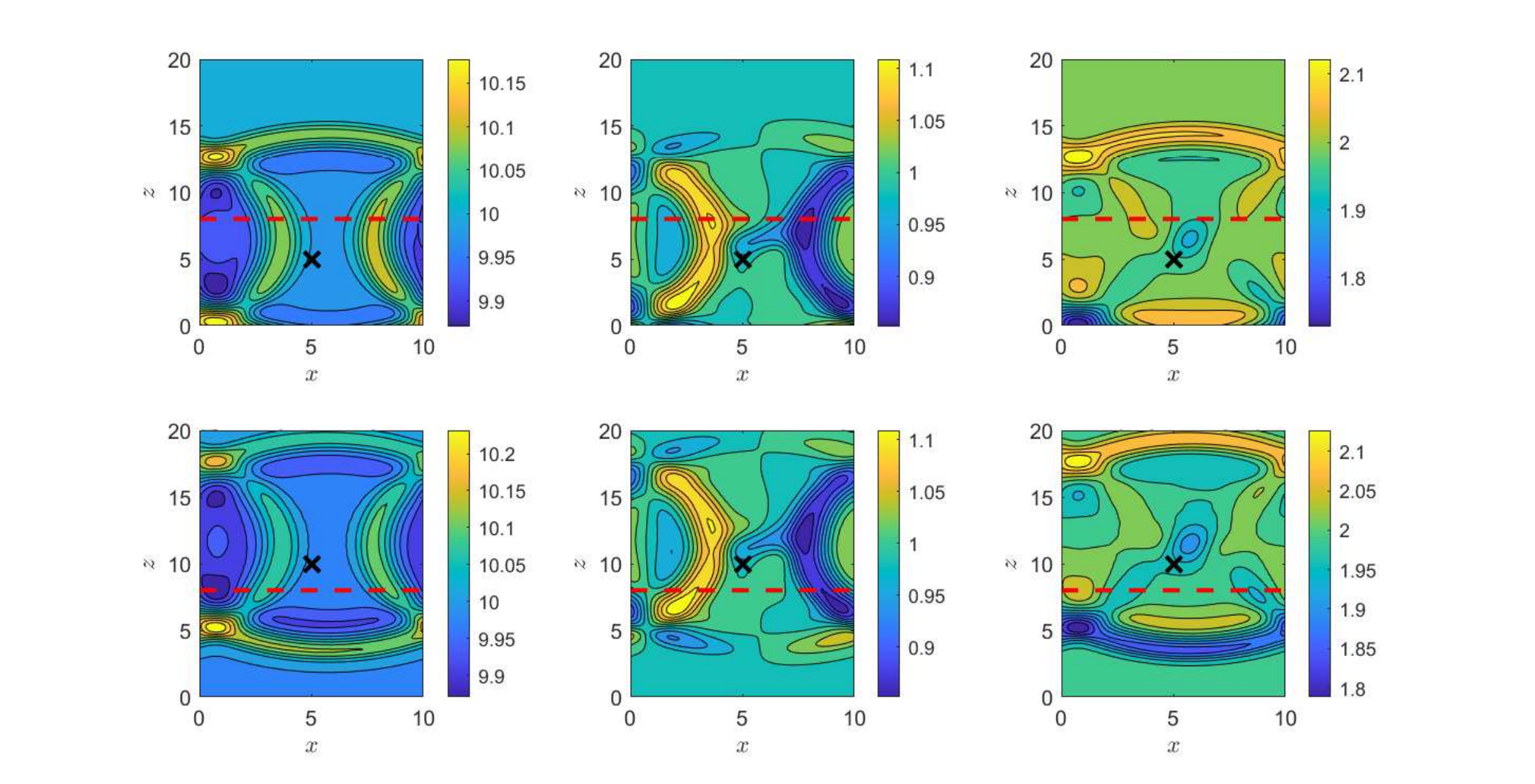}
    \caption{ Non-linear two-dimensional shallow water equations, Gaussian initial data. Contour plot (contour interval $0.05$) of the numerical solution of the XDG scheme at the final time $T=0.75\,s$ in the outgoing case (above) and incoming case (below) for Gaussian initial data. Depth $h(x,z,T)$ (left), horizontal velocity $u(x,z,T)$ (center) and vertical velocity $v(x,z,T)$ (right). $M=50$. The last Laguerre mode is placed at $44.52\,m$, only the region $[0,10\,m]\times[0,20\,m]$ is shown here. The black cross represents the center of the initial profile, while the red dashed line denotes the finite/semi-infinite interface.}
    \label{fig:swe_out_and_inc}
\end{figure}

The XDG scheme simulates wave motion at the correct speed both in the outgoing case and in the incoming case (Figure \ref{fig:swe_out_and_inc} for $M=50$). Relative errors compared to the reference single-domain DG solution are mostly far below one percent, and at most a few percent (Tables \ref{tab:swe_outgoing} and \ref{tab:swe_incoming}), thus further validating the coupling strategy in the case of hyperbolic systems.

\begin{table}[ht]
    \centering
    
    \begin{tabular}{cccccccc}
    \toprule
        $M$ & $\beta$ & $\mathcal{E}_2^{rel}(h)$ & $\mathcal{E}_\infty^{rel}(h)$ & $\mathcal{E}_2^{rel}(u)$ & $\mathcal{E}_\infty^{rel}(u)$ & $\mathcal{E}_2^{rel}(v)$ & $\mathcal{E}_\infty^{rel}(v)$ \\ \midrule
        10 & 1 & 5.24e-05 & 3.90e-04 & 1.01e-03 & 1.33e-02 & 2.36e-04 & 1.73e-03 \\
        30 & 3 & 7.02e-07 & 8.62e-06 & 4.58e-05 & 5.92e-04 & 3.82e-06 & 3.93e-05 \\
        50 & 5 & 2.08e-07 & 1.44e-06 & 1.11e-05 & 1.04e-04 & 9.63e-07 & 7.22e-06 \\
        \bottomrule
    \end{tabular}
    \caption{ Shallow water equations with Gaussian initial data, outgoing case. $L_x=10\,m$, $L_z=8\,m$, $N_x=50$, $N_z=40$, $p_x=p_z=1$, $T=0.75\,s$, $N_t=750$, $H=10\,m$, $U=1\,m/s$, $V=2\,m/s$. Relative errors in the finite region with respect to a single-domain DG discretization on $L_z^\prime=3L_z$ with $N_z^\prime=3N_z$. See also Figure \ref{fig:swe_out_and_inc}, top panel.}
    \label{tab:swe_outgoing}
\end{table}

\begin{table}[ht]
    \centering
    \begin{tabular}{cccccccc}\toprule
        $M$ & $\beta$ & $\mathcal{E}_2^{rel}(h)$ & $\mathcal{E}_\infty^{rel}(h)$ & $\mathcal{E}_2^{rel}(u)$ & $\mathcal{E}_\infty^{rel}(u)$ & $\mathcal{E}_2^{rel}(v)$ & $\mathcal{E}_\infty^{rel}(v)$ \\ \midrule
        10 & 1 & 1.60e-03 & 5.07e-03 & 6.02e-03 & 2.98e-02 & 7.45e-03 & 2.03e-02 \\ 
        30 & 3 & 1.32e-05 & 4.52e-05 & 2.16e-04 & 2.95e-03 & 6.29e-05 & 2.71e-04 \\
        50 & 5 & 1.24e-05 & 4.30e-05 & 4.69e-05 & 3.43e-04 & 5.88e-05 & 2.09e-04 \\
        \bottomrule
    \end{tabular}
    \caption{ Shallow water equations with Gaussian initial data, incoming case. $L_x=10\,m$, $L_z=8\,m$, $N_x=50$, $N_z=40$, $p_x=p_z=1$, $T=0.75\,s$, $N_t=750$, $H=10\,m$, $U=1\,m/s$, $V=2\,m/s$. Relative errors in the finite region with respect to a single-domain DG discretization on $L_z^\prime=3L_z$ with $N_z^\prime=3N_z$. See also Figure \ref{fig:swe_out_and_inc}, bottom panel.}
    \label{tab:swe_incoming}
\end{table}

\subsection{Extended DG vs. single domain DG performance comparison in large semi-infinite region}\label{subsec: undamped_runs_si_domain}
Next, we run a performance comparison between the XDG scheme and a single-domain DG discretization in wave simulation in the semi-infinite region. In these efficiency comparison tests, for simplicity we consider the undamped linear advection-diffusion equation with constant coefficients, $\gamma=0$ and Gaussian initial profile \eqref{eq:gauss_initcond}. The exact solution in $\mathbb{R}^2$ at a generic time $t$ is

\begin{equation}
    q(x,z,t) = \frac{A}{\sqrt{1+\displaystyle\frac{4\mu_x}{\sigma_x^2}t}\sqrt{1+\displaystyle\frac{4\mu_z}{\sigma_z^2}t}}\text{exp}\left[-\left(\frac{x-x_0-u_xt}{\sqrt{\sigma_x^2+4\mu_xt}}\right)^2\right]\text{exp}\left[-\left(\frac{z-z_0-u_zt}{\sqrt{\sigma_z^2+4\mu_zt}}\right)^2\right].
\end{equation}
\vspace{3mm}

\begin{figure}[ht]
\begin{minipage}[]{0.45\textwidth}
\centering
\includegraphics{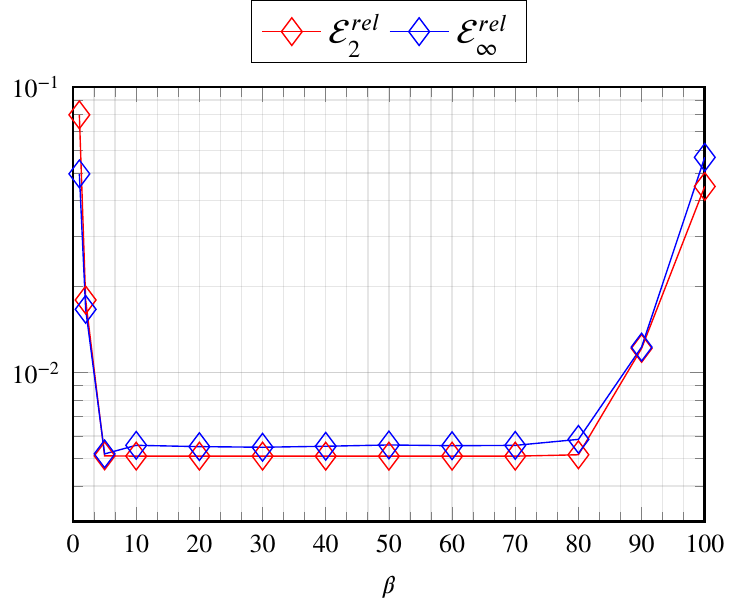}
\end{minipage}
\hfill
\begin{minipage}[]{0.55\textwidth}
\centering
    \begin{tabular}{c|c|c}
        $\beta$ & $z_M$ & $\delta$  \\\hline
        1 & 182.62 & 7.20e-02 \\
        2 & 91.31 & 3.60e-02 \\
        5 & 36.52 & 1.44e-02 \\
        10 & 18.26 & 7.20e-03 \\
        20 & 9.13 & 3.60e-03 \\
        30 & 6.09 & 2.40e-03 \\
        40 & 4.57 & 1.80e-03 \\
        50 & 3.65 & 1.44e-03 \\
        60 & 3.04 & 1.20e-03 \\
        70 & 2.61 & 1.03e-03 \\
        80 & 2.28 & 9.00e-04 \\
        90 & 2.03 & 8.00e-04 \\
        100 & 1.83 & 7.20e-04 \\
    \end{tabular}
\end{minipage}
\caption{Left: relative errors in the semi-infinite region as a function of $\beta$ for $M=50$. Right: position of the last Laguerre node ($z_M$) and distance between the first two Laguerre nodes ($\delta$) for the different choices of $\beta$.}\label{fig:err_vs_beta}
\end{figure}

\noindent The finite region is $[0,L_x]\times[0,L_z]$ with $L_x=0.5\,m$ and $L_z=1\,m$. The initial profile is placed inside the finite region ($x_0=0.25\,m$ and $z_0=0.5\,m$), with $\sigma_x=0.05\,m$ and $\sigma_z=0.15\,m$. By setting $u_x=0.5\,m/s$ and $u_z=2\,m$, the signal crosses the interface and ends up in the semi-infinite region at the final time $T=1\,s$. The viscosity coefficients in the $x$ and $z$ directions are $\mu_x=0.001\,m^2/s$ and $\mu_z=0.05\,m^2/s$. We compute relative errors with respect to the exact solution $q(x,z,T)$ by discretizing the semi-infinite part $[0,L_x]\times[L_z, \infty)$ in two different ways: (a) 2D XDG discretization using $M$ Laguerre nodes (Section \ref{sec:discrete_formulation}); 
(b) 2D DG discretization in $[0,L_x]\times[L_z,L_z+z_M]$ on a non-uniform grid in $z$, where the endpoints of the subintervals coincide with the Laguerre nodes in (a) and $z_M$ is the last Laguerre node.
The discretization parameters are $N_x=100$, $N_z=150$, $p_x=p_z=1$, $N_t=200$. \\ \indent In order to choose the value of $\beta$ for the simulations, we carry out an error analysis in the XDG scheme as $\beta$ varies, for a fixed value of $M$. Relative errors in the semi-infinite domain are computed at the final time with respect to the exact solution (Figure \ref{fig:err_vs_beta} left). We also report in Figure \ref{fig:err_vs_beta} (right) the extension of the semi-infinite domain and the distance between the first two Laguerre nodes for each choice of $\beta$. The analysis suggests that $\beta$ should be neither too large nor too small. In the former case, the spatial extension of the semi-infinite region is not large enough to accurately represent the solution beyond the interface, while in the latter case the first two nodes are too far apart, causing spurious signals due to the coupling. If $\beta$ is in those critical regions, the errors exhibit a noticeable increase. In addition, by fine-tuning the choice of $\beta$ the accuracy of the numerical solution can be improved by over one order of magnitude, with no effect on the total computational cost (Figure \ref{fig:err_vs_beta} left). On the basis of this analysis, for this set of tests we select the value $\beta=5$, which allows the representation of a large portion of the semi-infinite domain while keeping the errors as small as possible. 

We remark that techniques exist \cite{xia:2021} that further enhance the accuracy of Laguerre spectral discretizations by adjusting the extension of the semi-infinite domain through an adaptive choice of the $M$ and $\beta$ parameters. While such a feature would be of relatively little interest in the context of damped runs with a small number of Laguerre basis functions such as those considered in the previous sections, adaptive algorithms\cite{xia:2021} could be of great use where higher accuracy throughout the simulation might be required. The exploration of such techniques in the framework of the XDG model is left for future work.

\begin{figure}[htb!]
\centering
\includegraphics{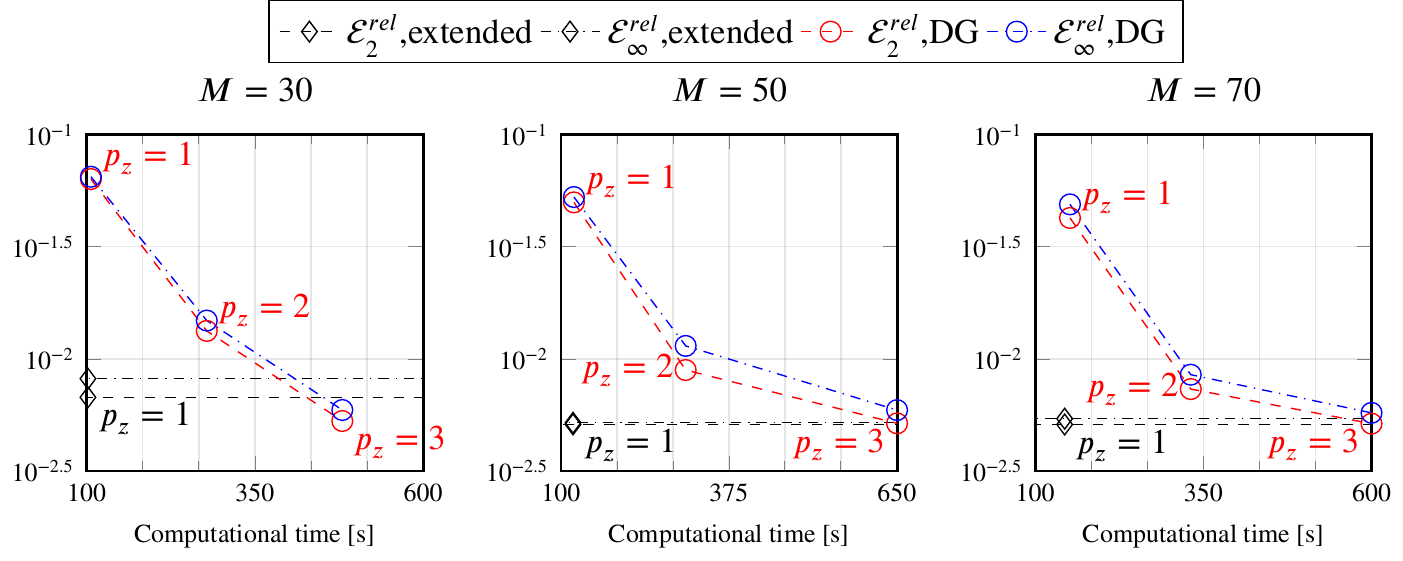}
\caption{Relative errors in the semi-infinite region vs. computational time, advection-diffusion of Gaussian initial data with $M=30$ (left), $M=50$ (center), $M=70$ (right) points in the semi-infinite domain, XDG scheme (black lines), single-domain DG scheme (red and blue lines). $L_x=0.5\,m$, $L_z=1\,m$, $N_x=100$, $p_x=1$, $N_z=150$, $\beta=5$, $T=1\,s$, $N_t=200$, $\mu_x=0.001\,m^2/s$, $\mu_z=0.05\,m^2/s$, $u_x=0.5\,m/s$, $u_z=2\,m/s$, $A=1\,m$, $x_0=0.25\,m$, $z_0=0.5\,m$, $\sigma_x=0.05\,m$, $\sigma_z=0.15\,m$.}\label{fig:err_infinite_region_beta5}
\end{figure}

With the choice $\beta=5$, errors obtained with the 2D XDG scheme are low even using a polynomial degree $p_z=1$ in the finite region (black lines in Figure \ref{fig:err_infinite_region_beta5}). By contrast, the 2D single-domain DG discretization on a non-uniform grid in $z$ achieves errors lower than one percent only using $p_z=3$ (red and blue lines in Figure \ref{fig:err_infinite_region_beta5}), at more than $4$ times the computational cost of the XDG model.

\begin{figure}[htb!]
\centering
\includegraphics{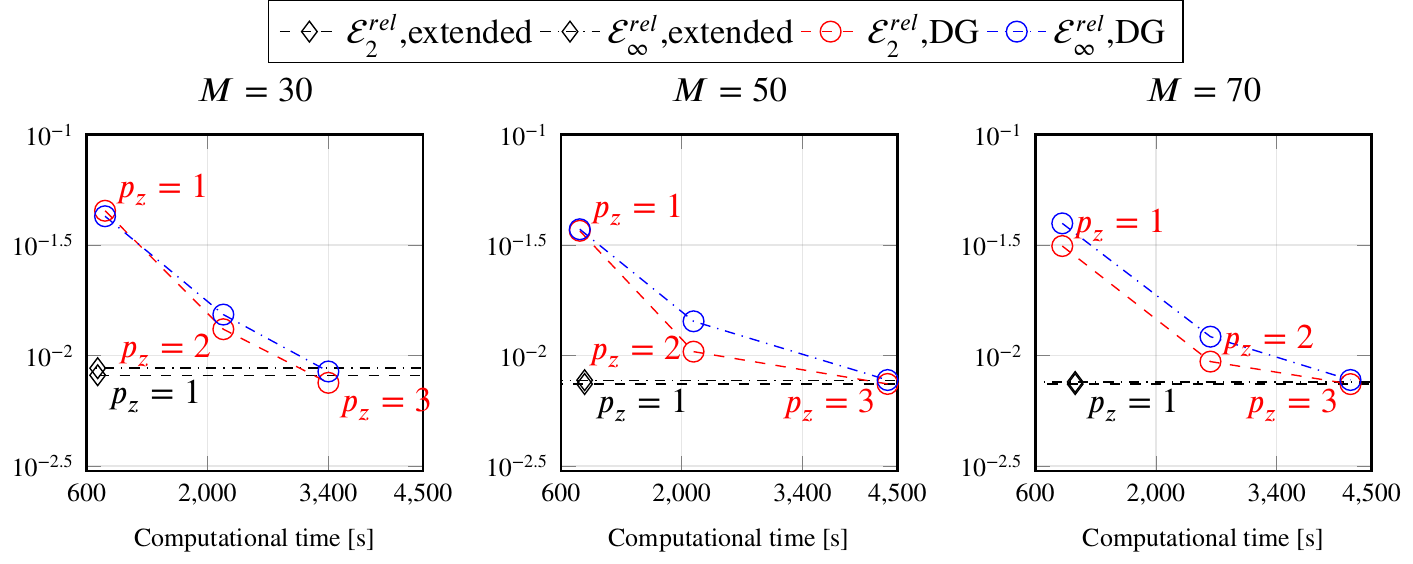}
\caption{Relative errors in the semi-infinite region vs. computational time, advection-diffusion of Gaussian initial data with $M=30$ (left), $M=50$ (center), $M=70$ (right) points in the semi-infinite domain, XDG scheme (black lines), single-domain DG scheme (red and blue lines). $L_x=0.5\,m$, $L_z=100\,m$, $N_x=100$, $p_x=1$, $N_z=150$, $\beta=0.05$, $T=7\,s$, $N_t=1400$, $\mu_x=5\times10^{-5}\,m^2/s$, $\mu_z=25\,m^2/s$, $u_x=1/7\,m/s$, $u_z=10\,m/s$, $A=1\,m$, $x_0=0.25\,m$, $z_0=80\,m$, $\sigma_x=0.05\,m$, $\sigma_z=15\,m$.}\label{fig:err_infinite_region_beta005}
\end{figure}

\newpage
Finally, we choose a smaller value of $\beta$ to simulate large spatial scales in the $z$ direction; more precisely, we set $\beta=0.05$ and we move the interface at $L_z=100\,m$. In this way, the last Laguerre node (i.e., the extension of the unbounded part of the XDG scheme) is placed at $z_M=633.7\,m$, $z_M=2120.9\,m$ and $z_M=3652.4\,m$ for $M=30,50,70$, respectively. With these parameter choices, the numerical solution obtained with the XDG scheme is at least as accurate as, and much less computationally costly than, the single-domain DG solution, with a significant efficiency gain (Figure \ref{fig:err_infinite_region_beta005} and Table \ref{tab:comptimes_comparison_beta005}). For $M=30$ points in the semi-infinite region, the single-domain DG scheme runs four times slower than the XDG scheme at given accuracy level. For $M=50$, the XDG scheme is five times more efficient.

\begin{table}[ht]
    \centering
    \begin{tabular}{ccccccc}
    \toprule
    & \multicolumn{2}{c}{$M=30$} & \multicolumn{2}{c}{$M=50$} & \multicolumn{2}{c}{$M=70$} \\ \cmidrule(lr){2-3}\cmidrule(lr){4-5}\cmidrule(lr){6-7}         & XDG & DG & XDG & DG & XDG & DG \\ \midrule
        $\mathcal{E}_2^{rel}$ & 8.14e-03 & 7.53e-03 & 7.43e-03 & 7.44e-03 & 7.42e-03 & 7.43e-03 \\
        $\mathcal{E}_\infty^{rel}$ & 8.76e-03 & 8.47e-03 & 7.73e-03 & 7.75e-03 & 7.59e-03 & 7.75e-03 \\
        $t$ & 727.8 & 3402.2 & 875.0 & 4386.8 & 1064.2 & 4255.4 \\
        $S$  & 4.67 & & 5.01 & & 4.00 & \\
        \bottomrule
    \end{tabular}
    \caption{ Comparison between the best relative errors obtained with the XDG scheme (XDG) and the single-domain DG scheme (DG) in the case $\beta=0.05$ and corresponding computational times $t$ in seconds. Polynomial degree $p_z=1$ is used for the XDG and $p_z=3$ for the single-domain DG discretization. The efficiency gain $S$ is defined as the ratio between the computational time of the single-domain DG run and that of the XDG run.}\label{tab:comptimes_comparison_beta005}
\end{table}

\newpage
\subsection{Efficiency of the XDG scheme in absorbing layer tests}
We now implement an absorbing layer in the semi-infinite region by setting $\gamma$ in \eqref{eq:adv_eq_coupled} as the sigmoid defined by
\begin{equation}
    \gamma(x,z) = \frac{\Delta\gamma}{1+\exp\left(\frac{\alpha L_0-z+L}{\sigma_D}\right)}.
\end{equation}
Here $\Delta\gamma$ is the sigmoid amplitude, $\alpha\in[0,1]$ the position of the sigmoid inside the absorbing layer, $L_0$ the spatial extension of the semi-infinite region and $\sigma_D$ the sigmoid steepness.

We consider the 2D nonlinear shallow water equations \eqref{eq:nl_swe} and impose a Dirichlet boundary condition sinusoidal in time to simulate a train of waves originating at the lower boundary $L_z=0$\cite{benacchio:2013,benacchio:2019}:
\begin{equation}
    q(x,0,t)=A\sin\left(\frac{2\kappa\pi t}{T}\right)\exp\left[-\left(\frac{x-x_0}{\sigma_0}\right)^2\right].
\end{equation}
The parameters in this experiment are $L_x=12\,m$, $L_z=2\,m$, $p_x=p_z=1$, $T=10\,s$, $N_t=10^4$. The sigmoid parameters are $\Delta\gamma=50\,m$, $\alpha=0.1$, $\sigma_D=30\,m$, while we set $A=1\,m$ in the boundary condition. The initial condition of the problem is $h(x,z,0)=10\,m$, $u(x,z,0)=v(x,z,0)=0$. The wave train crosses the finite region and is damped by the absorbing layer when it crosses the interface located at $L_z=2\,m$. We test the XDG strategy for different choices of $\kappa$, $N_z$, $N_x$ and $M$. In particular, we choose $(\kappa,N_z,N_x)=(100,20,120)$ or $(\kappa,N_z,N_x)=(200,40,240)$ and $M\in\{5,10,15\}$. For a fixed value of $\kappa$, $N_z$ and $N_x$, the value of the scaling parameter $\beta$ is chosen so that the distance between the first two Laguerre modes is approximately constant as $M$ varies. 

\begin{figure}[h]
    \centering
    \includegraphics[width=.925\textwidth]{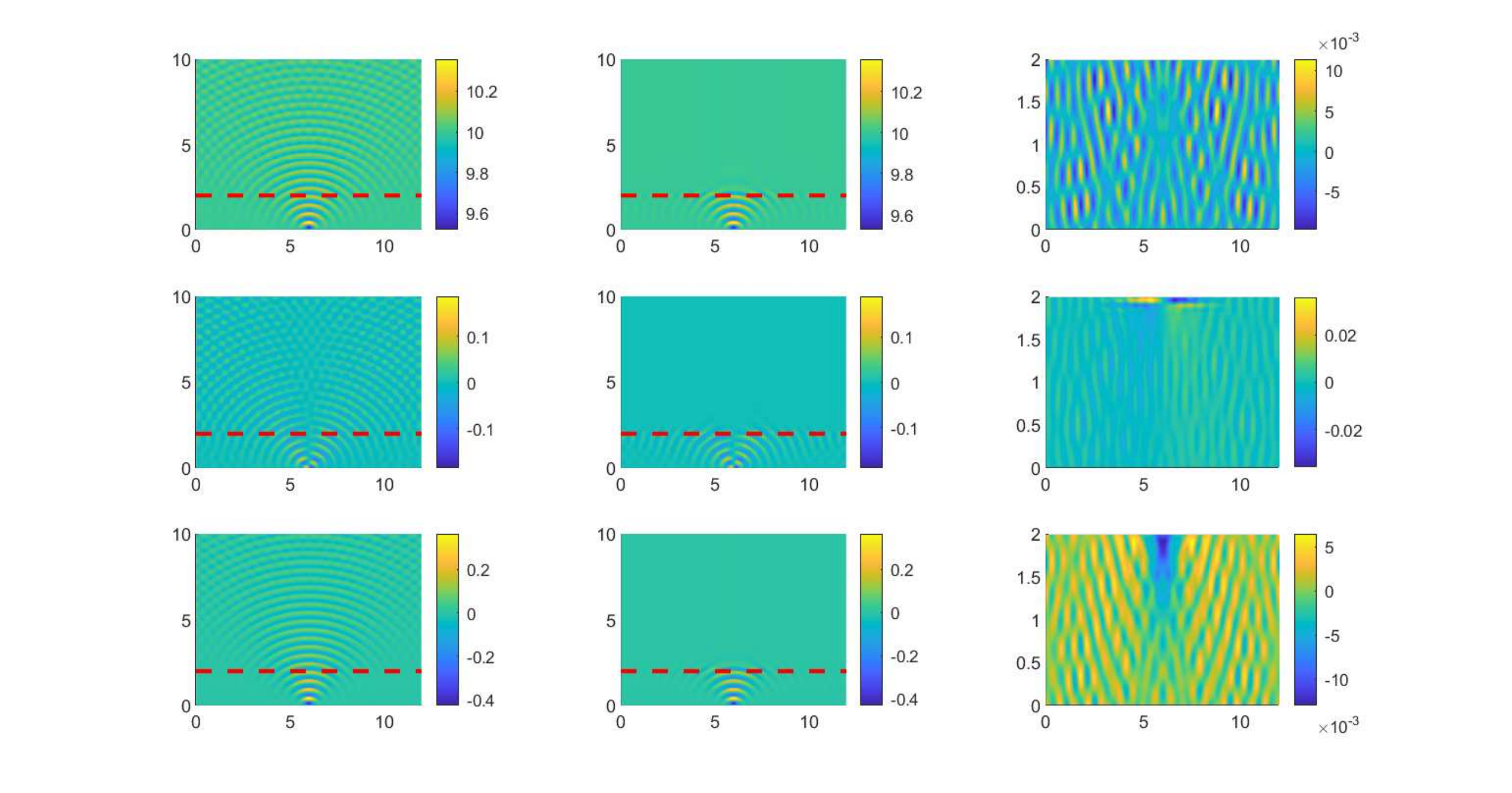}
    \caption{ Non-linear shallow water equations, wave train test. Single-domain DG solution using a uniform grid on $[0,12\,m]\times[0,10\,m]$ (left), solution of the XDG model with an absorbing layer in the semi-infinite region (center) and relative errors in the finite region $[0,12\,m]\times[0,2\,m]$ (right). $h(x,z,T)$ (above) $u(x,z,T)$ (center) and $v(x,z,T)$ (below). The red dashed line denotes the interface between finite and semi-infinite regions. $\kappa=200$, $M=5$, $\beta=8$.}
    \label{fig:wavetrain_plot_k200}
\end{figure}

\begin{figure}
    \centering
    \includegraphics[width=.925\textwidth]{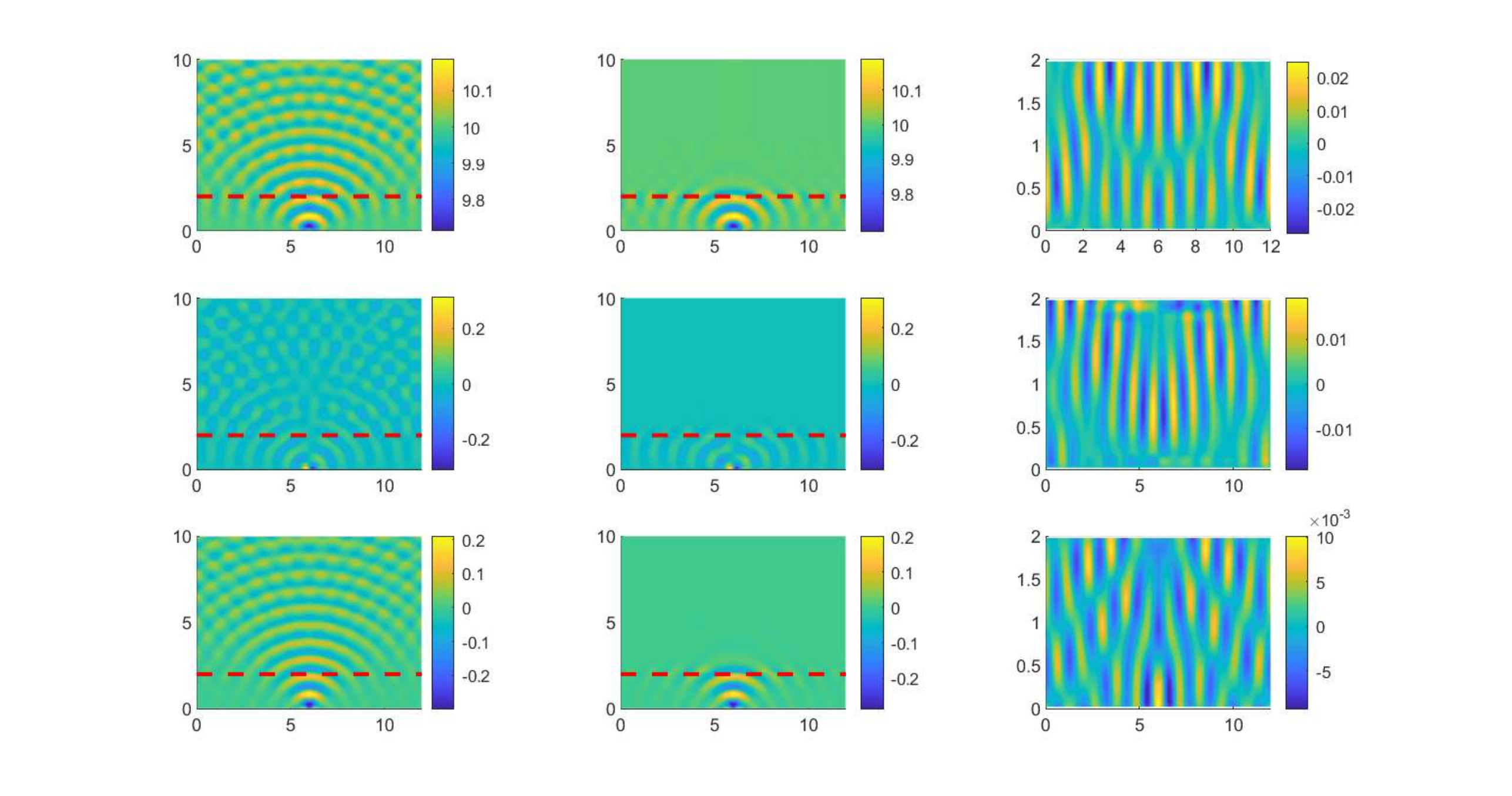}
    \caption{ Non-linear shallow water equations, wave train test. Single-domain DG solution using a uniform grid on $[0,12\,m]\times[0,10\,m]$ (left), solution of the XDG model with an absorbing layer in the semi-infinite region (center) and relative errors in the finite region $[0,12\,m]\times[0,2\,m]$ (right). $h(x,z,T)$ (above) $u(x,z,T)$ (center) and $v(x,z,T)$ (below). The red dashed line denotes the interface between finite and semi-infinite regions. $\kappa=100$, $M=5$, $\beta=4$.}
    \label{fig:wavetrain_plot_k100}
\end{figure}

First, to assess the efficiency of the absorbing layer, we compare the numerical solution of the 2D XDG scheme with a reference single-domain DG solution on $[0,L_x]\times[0,5L_z]$ using a \emph{uniform grid} with the same grid spacing $\Delta z$ as in the finite portion of the XDG scheme, and we compute errors at the final time. The absorbing layer, active in the semi-infinite part of the XDG setup, is able to damp outgoing perturbations (central panels of  Figures \ref{fig:wavetrain_plot_k200} and \ref{fig:wavetrain_plot_k100}). Spurious reflections in the finite region due to the use of two different sets of basis functions in the 2D XDG scheme result in relative errors in the finite region of a few percent at most (Table \ref{tab:wavetrain_err} and right panels of Figures \ref{fig:wavetrain_plot_k200} and \ref{fig:wavetrain_plot_k100}). Moreover, by suitably tuning the scaling parameter $\beta$, errors can be kept under control even with a small number of modes in the unbounded region. 

\begin{table}[ht]
    \centering
    \begin{tabular}{ccccccccccc}
    \toprule
        $\kappa$ & $N_z$ & $N_x$ & $M$ & $\beta$ & $\mathcal{E}_2^{rel}(h)$ & $\mathcal{E}_\infty^{rel}(h)$ &  $\mathcal{E}_2^{rel}(u)$ & $\mathcal{E}_\infty^{rel}(u)$ &  $\mathcal{E}_2^{rel}(v)$ & $\mathcal{E}_\infty^{rel}(v)$ \\ \midrule
          \multirow{3}{*}{200} & \multirow{3}{*}{40} & \multirow{3}{*}{240} & 15 & 3 & 1.17e-02 & 7.92e-03 & 1.80e-02 & 4.28e-02 & 1.02e-02 & 1.26e-02 \\
           & & & 10 & 4 & 1.24e-02 & 9.06e-03 & 1.95e-02 & 4.47e-02 & 1.14e-02 & 1.29e-02 \\
           & & & 5 & 8 & 1.56e-02 & 1.14e-02 & 1.89e-02 & 3.60e-02 & 1.18e-02 & 1.30e-02 \\ \hline
          \multirow{3}{*}{100} & \multirow{3}{*}{20} & \multirow{3}{*}{120} & 15 & 1.5 & 3.26e-02 & 2.16e-02 & 2.59e-02 & 1.78e-02 & 1.27e-02 & 8.00e-03 \\
           & & & 10 & 2 & 3.51e-02 & 2.33e-02 & 2.83e-02 & 1.93e-02 & 1.36e-02 & 8.64e-03 \\
           & & & 5 & 4 & 4.63e-02 & 2.77e-02 & 3.32e-02 & 1.91e-02 & 1.64e-02 & 1.01e-02
        \end{tabular}
    \caption{ Non-linear shallow water equations with an absorbing layer in the semi-infinite region, wave train test, XDG vs. uniform-grid DG. Relative errors of the XDG model at the final time with respect to a single-domain DG discretization in $[0,12\,m]\times[0,10\,m]$ with $N_z^\prime=5N_z$. (The grid spacing $\Delta z$ in the single-domain DG discretization is the same as that used in the finite region of the XDG scheme.)}
    \label{tab:wavetrain_err}
\end{table}

Second, to gauge the damping performance of the XDG scheme compared to a standard approach, we analyze the results of the single-domain DG discretization on a uniform grid with respect to a single-domain DG discretization on a \emph{non-uniform} grid in the $z$-direction, where the endpoints of the intervals coincide with the Laguerre nodes employed in the semi-infinite region for the XDG discretization. The relative errors in the finite region $[0,12\,m]\times[0,2\,m]$ (Table \ref{tab:wavetrain_err_DGnonunif}) are $O(10^{-2})$ or less, and comparable with those in Table \ref{tab:wavetrain_err}, showing that no accuracy penalty is incurred by using the XDG scheme as an absorbing layer. 

Although the final time for the tests in this section was kept short because of limited computational resources, additional numerical tests not shown here confirmed that the errors actually level off at the values shown in Tables \ref{tab:wavetrain_err} and \ref{tab:wavetrain_err_DGnonunif} when running until later times. In addition, the values did not vary significantly when refining the spatial or temporal spacing, thereby confirming that the errors are chiefly due to the presence of the absorbing layer and the choice of related parameters. 

While the two-dimensional results confirm the viability of the XDG scheme \cite{benacchio:2013,benacchio:2019,vismara:2022}, we expect the advantages of the XDG scheme to be more evident in problems with more complex dynamics. An in-depth exploration of the properties of the XDG scheme in damping internal waves arising as solutions of fluid flow equations in stratified media and non-flat bottom boundaries is left for future work.

\begin{table}[ht]
    \centering
    \begin{tabular}{ccccccccccc}
    \toprule
        $\kappa$ & $N_z$ & $N_x$ & $M$ & $\beta$ & $\mathcal{E}_2^{rel}(h)$ & $\mathcal{E}_\infty^{rel}(h)$ &  $\mathcal{E}_2^{rel}(u)$ & $\mathcal{E}_\infty^{rel}(u)$ &  $\mathcal{E}_2^{rel}(v)$ & $\mathcal{E}_\infty^{rel}(v)$ \\ \midrule
          \multirow{3}{*}{200} & \multirow{3}{*}{40} & \multirow{3}{*}{240} & 15 & 3 & 1.58e-02 & 9.95e-03 & 2.48e-02 & 3.08e-02 & 1.98e-02 & 2.67e-02 \\
           & & & 10 & 4 & 1.73e-02 & 1.03e-02 & 2.56e-02 & 3.13e-02 & 2.11e-02 & 2.87e-02 \\
           & & & 5 & 8 & 1.74e-02 & 1.16e-02 & 2.34e-02 & 2.95e-02 & 1.95e-02 & 2.51e-02 \\ \hline
          \multirow{3}{*}{100} & \multirow{3}{*}{20} & \multirow{3}{*}{120} & 15 & 1.5 & 2.97e-02 & 1.68e-02 & 2.37e-02 & 1.87e-02 & 1.86e-02 & 1.17e-02 \\
           & & & 10 & 2 & 3.98e-02 & 2.27e-02 & 2.23e-02 & 1.80e-02 & 1.95e-02 & 1.07e-02 \\
           & & & 5 & 4 & 4.00e-02 & 2.51e-02 & 2.50e-02 & 1.48e-02 & 1.97e-02 & 1.20e-02
        \end{tabular}
        \caption{Non-linear shallow water equations with an absorbing layer in the semi-infinite region, wave train test, non uniform-grid vs. uniform-grid DG. Relative errors of the single-domain DG model on a non-uniform grid  at the final time with respect to a single-domain DG discretization in $[0,12\,m]\times[0,10\,m]$ with $N_z^\prime=5N_z$.}
        \label{tab:wavetrain_err_DGnonunif}
\end{table}

\newpage
\section{Conclusion}\label{sec:disc_conc}

We proposed an extended discontinuous Galerkin numerical scheme for the discretization of hyperbolic and parabolic problems on two-dimensional semi-infinite domains. The method deals with the unbounded computational domain by splitting it into a finite region, where Legendre polynomials are used as basis and test functions, and an unbounded region, discretized using scaled Laguerre functions. The two subdomains are seamlessly coupled by means of numerical fluxes at the interface. For several integral terms within the discrete formulation, the expressions for the integrals was found to be formally identical except for replacing the finite subdomain grid spacing with the inverse Laguerre scaling parameter. Formal accuracy of the algorithm was retrieved in  convergence experiments, which confirmed the expected theoretical behaviour. No accuracy degradation issues were encountered linked to the tensor-product structure of the multi-dimensional discretization, obviating the need for approaches such as the hyperbolic cross approximation developed in \cite{shen:2010}.

Coupling validation tests were carried out on the 2D linear advection-diffusion equation, the 2D Burgers' equation, and the 2D non-linear shallow water equations using Gaussian initial data. Comparing the XDG solution with a standard single-domain DG discretization, we showed that spurious oscillations due to the coupling are negligible. Numerical experiments showed that the interface between the two regions is transparent to both outgoing and incoming signals, with minimal spurious reflections into the bounded domain.

Numerical evidence shows that the choice of the scaling parameter $\beta$ affects the accuracy of the solution in the semi-infinite part of the domain. This is an important feature of the method because the optimal selection of $\beta$ would allow one to achieve smaller errors with the same number of spectral modes $M$. While this work only reports the best results obtained by tuning $\beta$, a rigorous criterion for its selection remains an important open problem and shall be the subject of future work.

In addition, by suitable tuning of the scaling parameter of Laguerre basis functions, the XDG scheme proved to be a highly efficient tool to accurately represent dynamics over arbitrarily large spatial scales. Achieving comparable accuracy in the semi-infinite region with a single-domain DG discretization required higher polynomial degrees and more than four times as much computational time in performance comparison tests.

By means of a sigmoidal reaction term, an absorbing layer was then implemented in the semi-infinite region. Exponential convergence provided by the Laguerre basis functions enabled the scheme to efficiently damp outgoing wave trains with a small number of modes. Reflections into the finite domain measured as errors with respect to a single-domain DG discretization on a uniform grid were of comparable amplitude to those obtained by using a standard DG approach on a non-uniform grid with nodes coinciding with the Laguerre nodes in the semi-infinite region.

Numerical evidence suggests that the proposed XDG approach can be a competitive tool to discretize transient dynamics efficiently over arbitrarily large computational domains. Inspired by these results, several research directions can be explored in future work. First, a natural extension would be the use of polar coordinates to simulate semi-infinite circular domains. In that setup, more complex nonlinear magnetohydrodynamics models could be considered in order to gauge the potential of the XDG scheme to model phenomena spanning very large length scales such as coronal mass ejections. Second, with the addition of bottom topography, the scheme could be used to discretize compressible Euler equations first in vertical slice two-dimensional domains and then in full-fledged three-dimensional implementations. Established idealized orography benchmarks \cite{bonaventura:2000,melvin:2019} could test the method's effectiveness as an efficient damping layer. If the techniques proposed in our work achieve the same accuracy at a lower computational cost than currently employed approaches, related savings can be reinvested in the region of interest to achieve e.g. higher accuracy or longer lead times in forecasts.

Finally, increasing resolution in the unbounded region would enable comparison with extensions of weather forecast models to the middle and upper atmosphere and, in perspective, use of the method in space weather models\cite{borchert:2019,jackson:2019,klemp:2021}. Future applications of the method could also benefit from a recently proposed adaptive capability \cite{chou:2022,xia:2021} whereby the scaling parameter can vary within simulations, thus enabling the scheme to simulate multiple phenomena developing over very different spatial and temporal scales.

\section*{Authorship contribution statement}

\textbf{F. Vismara:} conceptualization, investigation, methodology, software, validation, visualization, writing - original draft, writing - review \& editing, data curation. \textbf{T. Benacchio:} conceptualization, investigation, methodology, writing - original draft, writing - review \& editing, supervision, data curation.

\section*{Acknowledgements}
The authors would like to thank Prof. Luca Bonaventura for useful feedback on early drafts of the paper.

\section*{Data availability}

The datasets generated during the current study are available from the corresponding author on reasonable request.

\appendix 

\section{Computation of the integrals in the discrete XDG formulation}
\label{app: discrete_form}

We analyze in detail all terms appearing in \eqref{eq:weak_form_final_disc_coupled}. For simplicity we consider a uniform grid $\Delta x_{m_x}=\Delta x, \forall m_x$, $\Delta z_{m_z}=\Delta z, \forall m_z$. Using the fact that
\begin{equation}\label{eq:laglegformulas_coupled}
    \int_{x_{m_x-1}}^{x_{m_x}}\phi^{m_x}_j(x)\phi^{m_x}_{j^\prime}(x)dx=\Delta x\delta_{jj^\prime},\qquad\int_{L_z}^{+\infty}\phi_i^\infty(z)\phi_{i^\prime}^\infty(z) dz=\frac{1}{\beta}\delta_{ii^\prime},
\end{equation}
the first integral can easily be written as
\begin{equation}\label{eq:1st_coupled}
    \frac{d}{dt}\int_{K_{m_x,m_z}}\sum_{j}\sum_{i}q_{m_x,m_z}^{(j,i)}(t)\phi_j^{m_x}(x)\phi_i^{m_z}(z)\phi^{m_x^\prime}_{j^\prime}(x)\phi_{i^\prime}^{m_z^\prime}(z)dx\,dz = \begin{cases} \Delta x\Delta z\displaystyle\frac{dq_{m_x^\prime,m_z^\prime}^{(j^\prime,i^\prime)}}{dt} & m_z^\prime=1,\dots,N_z \\ \displaystyle\frac{\Delta x}{\beta}\displaystyle\frac{dq_{m_x^\prime,m_z^\prime}^{(j^\prime,i^\prime)}}{dt} & m_z^\prime=\infty \end{cases}.
\end{equation}

\noindent The fifth integral is
\begin{equation}
    \int_\Omega\bm{\mu}\nabla q_h\cdot\nabla\varphi_h\,dx\,dz = \int_\Omega\left(\mu_x\frac{\partial q_h}{\partial x}\frac{\partial\varphi_h}{\partial x}+\mu_z\frac{\partial q_h}{\partial z}\frac{\partial\varphi_h}{\partial z}\right)\,dx\,dz
\end{equation}
(notice that, under the assumption that $\bm{\mu}$ is diagonal, $\bm{\mu}\nabla=[\mu_x\partial/\partial x,\mu_z\partial/\partial z]$). In $\Omega_{DG-DG}$ ($m_z=1,\dots,N_z$) this gives
\begin{align}\label{eq:5th_coupled_DGDG}
\begin{split}
    &\mu_x\sum_{j=0}^{p_x}\sum_{i=0}^{p_z}q_{m_x,m_z}^{(j,i)}\int_{K_{m_x,m_z}}(\phi^{m_x}_j)^\prime\phi^{m_z}_i(\phi^{m_x}_{j^\prime})^\prime\phi^{m_z}_{i^\prime}+\mu_z\sum_{j=0}^{p_x}\sum_{i=0}^{p_z}q_{m_x,m_z}^{(j,i)}\int_{K_{m_x,m_z}}\phi^{m_x}_j(\phi^{m_z}_i)^\prime\phi^{m_x}_{j^\prime}(\phi^{m_z}_{i^\prime})^\prime=\\&=\mu_x\Delta z\sum_{j=0}^{p_x}q_{m_x,m_z}^{(j,i^\prime)}\int_{x_{m_x-1}}^{x_{m_x}}(\phi^{m_x}_j)^\prime(x)(\phi^{m_x}_{j^\prime})^\prime(x)\,dx+\mu_z\Delta x\sum_{i=0}^{p_z}q_{m_x,m_z}^{(j^\prime,i)}\int_{z_{m_z-1}}^{z_{m_z}}(\phi^{m_z}_i)^\prime(z)(\phi^{m_z}_{i^\prime})^\prime(z)\,dz,
\end{split}
\end{align}
while in $\Omega_{DG-LAG}$ ($m_z=\infty$) we have 
\begin{align}\label{eq:5th_coupled_DGLAG}
\begin{split}
    &\mu_x\sum_{j=0}^{p_x}\sum_{i=0}^{M}q_{m_x,\infty}^{(j,i)}\int_{K_{m_x,\infty}}(\phi^{m_x}_j)^\prime\phi^{\infty}_i(\phi^{m_x}_{j^\prime})^\prime\phi^{\infty}_{i^\prime}+\mu_z\sum_{j=0}^{p_x}\sum_{i=0}^{M}q_{m_x,\infty}^{(j,i)}\int_{K_{m_x,\infty}}\phi^{m_x}_j(\phi^{\infty}_i)^\prime\phi^{m_x}_{j^\prime}(\phi^{\infty}_{i^\prime})^\prime=\\&=\frac{\mu_x}{\beta}\sum_{j=0}^{p_x}q_{m_x,\infty}^{(j,i^\prime)}\int_{x_{m_x-1}}^{x_{m_x}}(\phi^{m_x}_j)^\prime(x)(\phi^{m_x}_{j^\prime})^\prime(x)\,dx+\mu_z\Delta x\sum_{i=0}^{M}q_{m_x,\infty}^{(j^\prime,i)}\int_{L_z}^{\infty}(\phi^{\infty}_i)^\prime(z)(\phi^{\infty}_{i^\prime})^\prime(z)\,dz,
\end{split}
\end{align}
i.e., the same expression as \eqref{eq:5th_coupled_DGDG} with $1/\beta$ replacing $\Delta z$, $\sum_{i=0}^M$ replacing $\sum_{i=0}^{p_z}$, $q_{m_x, \infty}$ replacing $q_{m_x,m_z}$, $K_{m_x, \infty}$ replacing $K_{m_x,m_z}$, and $\int_{L_z}^\infty$ replacing  $\int_{z_{m_z-1}}^{z_{m_z}}$.
Notice that, thanks to the definition of $\phi^\infty_i$,
\begin{equation}
    \int_{L_z}^{\infty}(\phi^{\infty}_i)^\prime(z)(\phi^{\infty}_{i^\prime})^\prime(z)\,dz=\int_0^\infty(\hat{\mathscr{L}}_i^\beta)^\prime(z)(\hat{\mathscr{L}}_{i^\prime}^\beta)^\prime(z)\,dz.
\end{equation}

\noindent In a similar fashion, the seventh integral is
\begin{equation}
    -\int_\Omega \mathbf{F}(q_h)\cdot\nabla\varphi_h\,dx\,dz=-\int_\Omega\left(F_1(q_h)\frac{\partial\varphi_h}{\partial x} + F_2(q_h)\frac{\partial\varphi_h}{\partial z}\right)\,dx\,dz,
\end{equation}
which locally reads, for all $m_x$ and $m_z$,
\begin{align}\label{eq:7th_coupled}
    \begin{split}
    -\int_{K_{m_x,m_z}}F_1\left(\sum_{j=0}^{p_x}\sum_{i}q_{m_x,m_z}^{(j,i)}\phi^{m_x}_j\phi^{m_z}_i\right)(\phi^{m_x}_{j^\prime})^\prime\phi^{m_z}_{i^\prime}-\int_{K_{m_x,m_z}}F_2\left(\sum_{j=0}^{p_x}\sum_iq_{m_x,m_z}^{(j,i)}\phi^{m_x}_j\phi^{m_z}_i\right)\phi^{m_x}_{j^\prime}(\phi^{m_z}_{i^\prime})^\prime.
    \end{split}
\end{align}
The sum over $i$ runs from $0$ to $p_z$ if $m_z=1,\dots,N_z$, from $0$ to $M$ if $m_z=\infty$.

Next, we consider the edge integrals, i.e. the second, third, fourth and sixth integral in \eqref{eq:weak_form_final_disc_coupled}; to this end we introduce the quantities
\begin{gather}
    \phi_j^L=\phi_j^{m_x}(x_{m_x-1}^+)=\phi_j^{m_z}(z_{m_z-1}^+) \qquad \phi_j^R=\phi_j^{m_x}(x_{m_x}^-)=\phi_j^{m_z}(z_{m_z}^-) \\
    (\phi_j^\prime)^L=(\phi_j^{m_x})^\prime(x_{m_x-1}^+)=(\phi_j^{m_z})^\prime(z_{m_z-1}^+) \qquad (\phi_j^\prime)^R=(\phi_j^{m_x})^\prime(x_{m_x}^-)=(\phi_j^{m_z})^\prime(z_{m_z}^-)
\end{gather}
for $m_z=1,\dots,N_z$, which represent the value of basis functions and their derivatives at the left and right endpoint of each subintervals. Notice that, on a uniform grid, these quantities do not depend on $m_x$ and $m_z$.

We start from the second, third and fourth term in \eqref{eq:weak_form_final_disc_coupled}. For the sake of exposition, in the following we specialize the formulation by considering horizontal edges aligned with the $x$ direction and vertical edges aligned with the $z$ direction.
Let $e\in\Gamma_I^{DG-DG}$ be a \emph{vertical} internal edge in $\Omega_{DG-DG}$. Then, $e=\{x_{m_x}\}\times[z_{m_z-1},z_{m_z}]$ for $m_x=1,\dots,N_x-1$ and $m_z=1,\dots,N_z-1$, with $\mathbf{n}=[1,0]$ and $\lvert e\rvert=\Delta z$, and the corresponding terms are
\begin{equation}\label{}
    -\mu_x\int_{z_{m_z-1}}^{z_{m_z}}\left\{\left\{\frac{\partial q_h}{\partial x}\right\}\right\}[[\varphi_h]]\,dz + \mu_x\epsilon\int_{z_{m_z-1}}^{z_{m_z}}\left\{\left\{\frac{\partial \varphi_h}{\partial x}\right\}\right\}[[q_h]]\,dz + \frac{\sigma}{\Delta z}\int_{z_{m_z-1}}^{z_{m_z}}[[q_h]][[\varphi_h]]\,dz.
\end{equation}
Expanding $q_h$ and $\varphi_h$, after some derivations the three integrals above become
\begin{equation}\label{eq:diff_internal_DGDG_vertical_coupled}
    \sum_{j=0}^{p_x}\left(q_{m_x,m_z}^{(j,i^\prime)}m^{V,11}_{j,j^\prime}\delta_{m_x^\prime,m_x}+q_{m_x,m_z}^{(j,i^\prime)}m^{V,12}_{j,j^\prime}\delta_{m_x^\prime,m_x+1}+q_{m_x+1,m_z}^{(j,i^\prime)}m^{V,21}_{j,j^\prime}\delta_{m_x^\prime,m_x}+q_{m_x+1,m_z}^{(j,i^\prime)}m^{V,22}_{j,j^\prime}\delta_{m_x^\prime,m_x+1}\right).
\end{equation}

\noindent The local matrix $m^{V,11}$ (resp., $m^{V,22}$) corresponds to the interaction of the left neighboring element (resp., of the right neighboring element) with itself, while the local matrices $m^{V,12}$ and $m^{V,21}$ correspond to the interaction across the interface.

We obtain a similar result when $e\in\Gamma_I^{DG-DG}$ is a \emph{horizontal} internal edge in $\Omega_{DG-DG}$; in this case $e=[x_{m_x-1},x_{m_x}]\times\{z_{m_z}\}$ for $m_x=1,\dots,N_x-1$ and $m_z=1,\dots,N_z-1$, while $\mathbf{n}=[0,1]$ and $\lvert e\rvert=\Delta x$. Similar computations as above lead to
\begin{equation}\label{eq:diff_internal_DGDG_horizontal_coupled}
    \sum_{i=0}^{p_z}\left(q_{m_x,m_z}^{(j^\prime,i)}m^{H,11}_{i,i^\prime}\delta_{m_z^\prime,m_z}+q_{m_x,m_z}^{(j^\prime,i)}m^{H,12}_{i,i^\prime}\delta_{m_z^\prime,m_z+1}+q_{m_x,m_z+1}^{(j^\prime,i)}m^{H,21}_{i,i^\prime}\delta_{m_z^\prime,m_z}+q_{m_x,m_z+1}^{(j^\prime,i)}m^{H,22}_{i,i^\prime}\delta_{m_z^\prime,m_z+1}\right).
\end{equation}

\noindent We then consider $e\in\Gamma_I^{DG-LAG}$. Note that there are no horizontal edges here. Expanding the second, third and fourth integrals of \eqref{eq:weak_form_final_disc_coupled} for this case leads, after some derivations, to:
\begin{equation}\label{eq:diff_internal_DGLAG_coupled}
    \sum_{j=0}^{p_x}\left(q_{m_x,\infty}^{(j,i^\prime)}m^{V,11}_{j,j^\prime}\delta_{m_x^\prime,m_x}+q_{m_x,\infty}^{(j,i^\prime)}m^{V,12}_{j,j^\prime}\delta_{m_x^\prime,m_x+1}+q_{m_x+1,\infty}^{(j,i^\prime)}m^{V,21}_{j,j^\prime}\delta_{m_x^\prime,m_x}+q_{m_x+1,\infty}^{(j,i^\prime)}m^{V,22}_{j,j^\prime}\delta_{m_x^\prime,m_x+1}\right).
\end{equation}

\noindent Finally, we analyze the exchange of information at the interface by considering $e\in\Gamma_I^{INTERF}$. Since these edges are horizontal, we have
\begin{equation}
    -\mu_z\int_{x_{m_x-1}}^{x_{m_x}}\left\{\left\{\frac{\partial q_h}{\partial z}\right\}\right\}[[\varphi_h]]\,dx + \mu_z\epsilon\int_{x_{m_x-1}}^{x_{m_x}}\left\{\left\{\frac{\partial \varphi_h}{\partial z}\right\}\right\}[[q_h]]\,dx + \frac{\sigma}{\Delta x}\int_{x_{m_x-1}}^{x_{m_x}}[[q_h]][[\varphi_h]]\,dx.
\end{equation}

\noindent After some derivations, the above expression can be written compactly as
\begin{equation}\label{eq:diff_interface_DGLAG_coupled}
    \sum_{i=0}^{p_z}\left(q_{m_x,N_z}^{(j^\prime,i)}m^{C,11}_{i,i^\prime}\delta_{m_z^\prime,N_z}+q_{m_x,N_z}^{(j^\prime,i)}m^{C,12}_{i,i^\prime}\delta_{m_z^\prime,\infty}\right)+\sum_{i=0}^M\left(q_{m_x,\infty}^{(j^\prime,i)}m^{C,21}_{i,i^\prime}\delta_{m_z^\prime,N_z}+q_{m_x,\infty}^{(j^\prime,i)}m^{C,22}_{i,i^\prime}\delta_{m_z^\prime,\infty}\right)
\end{equation}
and the coupling between the two discretizations is described by the matrices $m^C$ (see Appendix \ref{sec: appendix_m}).

To conclude the analysis of the second, third and fourth terms of \eqref{eq:weak_form_final_disc_coupled} we consider Dirichlet edges $e\in\Gamma_D$, where we need to apply the definition of jumps and averages on boundary edges. If $e$ belongs to the left boundary $x=0$, then $\mathbf{n}=[-1,0]$ and we find
\begin{align}\label{eq:diff_dirichlet_coupled_left_DGDG}
\begin{split}
    &-\mu_x\int_{z_{m_z-1}}^{z_{m_z}}-\frac{\partial {q_h}_{|_e}}{\partial x}{\varphi_h}_{|_e}\,dz+\mu_x\epsilon\int_{z_{m_z-1}}^{z_{m_z}}-\frac{\partial{\varphi_h}_{|_e}}{\partial x}{q_h}_{|_e}\,dz+\frac{\sigma}{\Delta z}\int_{z_{m_z-1}}^{z_{m_z}}{q_h}_{|_e}{\varphi_h}_{|_e}\,dz=\\&=\sum_{j=0}^{p_x}\sum_{i=0}^{p_z}q_{1,m_z}^{(j,i)}\left(\mu_x\int_e(\phi_j^\prime)^L\phi_i^{m_z}\phi_{j^\prime}^L\phi_{i^\prime}^{m_z}-\mu_x\epsilon\int_e(\phi_{j^\prime}^\prime)^L\phi_{i^\prime}^{m_z}\phi_j^L\phi_i^{m_z}+\frac{\sigma}{\Delta z}\int_e\phi_j^L\phi_i^{m_z}\phi_{j^\prime}^L\phi_{i^\prime}^{m_z}\right)\delta_{m_x^\prime,1}=\\&=
    \sum_{j=0}^{p_x}q_{1,m_z}^{(j,i^\prime)}\left(\mu_x\Delta z(\phi_j^\prime)^L\phi_{j^\prime}^L-\mu_x\epsilon\Delta z\phi_j^L(\phi_{j^\prime}^\prime)^L+\sigma\phi_j^L\phi_{j^\prime}^L\right)\delta_{m_x^\prime,1}
\end{split}
\end{align}
for $m_z=1,\dots,N_z$ and a similar expression for $m_z=\infty$ with modifications as above. Analogous computations on the right boundary $x=L_x$, where $\mathbf{n}=[1,0]$, lead to
\begin{equation}\label{eq:diff_dirichlet_coupled_right_DGDG}
    \sum_{j=0}^{p_x}q_{N_x,m_z}^{(j,i^\prime)}\left(-\mu_x\Delta z(\phi_j^\prime)^R\phi_{j^\prime}^R+\mu_x\epsilon\Delta z\phi_j^R(\phi_{j^\prime}^\prime)^R+\sigma\phi_j^R\phi_{j^\prime}^R\right)\delta_{m_x^\prime,N_x}
\end{equation}
for $m_z=1,\dots,N_z$; the expression for $m_z=\infty$ is identical as long as one replaces $\Delta z$ with $1/\beta$. 

Finally, if $e$ belongs to the bottom boundary, then the integrals in \eqref{eq:weak_form_final_disc_coupled} only involve the DG-DG section and we obtain
\begin{equation}\label{eq:diff_dirichlet_coupled_bottom}
    \sum_{i=0}^{p_z}q_{m_x,1}^{(j^\prime,i)}\left(\mu_z\Delta x(\phi_i^\prime)^L\phi_{i^\prime}^L-\mu_z\epsilon\Delta x\phi_i^L(\phi_{i^\prime}^\prime)^L+\sigma\phi_i^L\phi_{i^\prime}^L\right)\delta_{m_z^\prime,1}.
\end{equation}

In summary, equations \eqref{eq:diff_internal_DGDG_vertical_coupled},\eqref{eq:diff_internal_DGDG_horizontal_coupled},\eqref{eq:diff_dirichlet_coupled_left_DGDG}\eqref{eq:diff_dirichlet_coupled_right_DGDG},\eqref{eq:diff_dirichlet_coupled_bottom} involve the DG-DG section, equations \eqref{eq:diff_internal_DGLAG_coupled} and \eqref{eq:diff_dirichlet_coupled_left_DGDG}-\eqref{eq:diff_dirichlet_coupled_right_DGDG} with $m_z=\infty$ and $\Delta z=1/\beta$ involve the DG-LAG section, and equation \eqref{eq:diff_interface_DGLAG_coupled} describes the coupling between the two sections at the interface $\Gamma_I^{INTERF}$. This concludes the analysis of terms two, three and four of equation \eqref{eq:weak_form_final_disc_coupled}.

To expand the sixth integral we need to introduce a suitable numerical flux at the internal edges $e$. We use the Rusanov flux, which for a general interface $e$ is defined as
\begin{equation}\label{eq:rusanov_coupled}
    \hat{q}({q_h}_{|_e},{\mathbf{n}}_{|_e})=\frac{1}{2}(\mathbf{F}^a(q_h^a) + \mathbf{F}^b(q_h^b))\cdot\mathbf{n}-\frac{1}{2}\nu_e(q_h^b-q_h^a),
\end{equation}
where $\mathbf{F}^a(q)=F(q,x^a,z)$, $\mathbf{F}^b(q)=F(q,x^b,z)$ if the interface is vertical and $\mathbf{F}^a(q)=F(q,x,z^a)$, $\mathbf{F}^b(q)=F(q,x,z^b)$ if it is horizontal, where $x^a:=x_e^-$, $x^b:=x_e^+$, $z^a:=z_e^-$, $z^b:=z_e^+$. Moreover, 
\begin{equation}
    \nu_e=\text{max}\left(\left\lvert\frac{d\mathbf{F}^a}{dq}(q_h^a)\cdot\mathbf{n}\right\rvert,\left\lvert\frac{d\mathbf{F}^b}{dq}(q_h^b)\cdot\mathbf{n}\right\rvert\right),
\end{equation}
where, for example, $\frac{d\mathbf{F}}{dq}=\mathbf{u}$ for the linear advection-diffusion equation and $\frac{d\mathbf{F}}{dq}(q)=[q,q]^\top$ for the Burgers' equation. Here we use the notation $\frac{d\mathbf{F}^a}{dq}(q)=\frac{d\mathbf{F}}{dq}(q,x^a,z)$ for vertical interfaces and similarly in the other cases.

We begin the analysis by considering internal edges. Notice that, when expanding the sum in the sixth term of \eqref{eq:weak_form_final_disc_coupled}, each internal edge is considered twice, once for each element to which it belongs. If $e$ is a \emph{vertical} internal edge, then we have the two expressions
\begin{equation}
    \int_e\left(\frac{1}{2}(F_1^a(q_h^a)+F_1^b(q_h^b))-\frac{1}{2}\nu_e(q_h^b-q_h^a)\right)\varphi_h^a\,dz \quad
    \int_e\left(\frac{1}{2}(-F_1^a(q_h^a)-F_1^b(q_h^b))-\frac{1}{2}\nu_e(q_h^a-q_h^b)\right)\varphi_h^b\,dz
\end{equation}

\noindent for $\mathbf{n}=[1,0]$ and $\mathbf{n}=[-1,0]$, respectively. Putting the two terms together we obtain
\begin{equation}\label{eq:adv_ver_interfaces}
    \int_e\left(\frac{F_1^a(q_h^a)+\nu_eq_h^a}{2}\varphi_h^a+\frac{-F_1^a(q_h^a)-\nu_eq_h^a}{2}\varphi_h^b+\frac{F_1^b(q_h^b)-\nu_eq_h^b}{2}\varphi_h^a+\frac{-F_1^b(q_h^b)+\nu_eq_h^b}{2}\varphi_h^b\right)\,dz
\end{equation}
and $\nu_e=\max(\lvert \frac{dF_1^a}{dq}(q_h^a)\rvert,\lvert \frac{dF_1^b}{dq}(q_h^b)\rvert)$. We can now insert the representation of $q_h$ and $\varphi_h$ and numerically compute the integral. 

If \emph{horizontal} internal edges $e$ are considered, the normals are $[0,-1]$ or $[0,1]$, so that the Rusanov flux yields
\begin{equation}\label{eq:adv_hor_interfaces}
    \int_e\left(\frac{F_2^a(q_h^a)+\nu_eq_h^a}{2}\varphi_h^a+\frac{-F_2^a(q_h^a)-\nu_eq_h^a}{2}\varphi_h^b+\frac{F_2^b(q_h^b)-\nu_eq_h^b}{2}\varphi_h^a+\frac{-F_2^b(q_h^b)+\nu_eq_h^b}{2}\varphi_h^b\right)\,dz,
\end{equation}
where now $\nu_e=\max(\lvert \frac{dF_2^a}{dq}(q_h^a)\rvert,\lvert \frac{dF_2^b}{dq}(q_h^b)\rvert)$ and  $e=[x_{m_x-1},x_{m_x}]\times\{z_{m_z}\}$, with $F_2^a(q_h^a)=F_2(q_h^a,x,z_{m_z}^-)$ and $F_2^b(q_h^b)=F_2(q_h^b,x,z_{m_z}^+)$, see Appendix \ref{sec: appendix_flux_int} for the expressions of unknowns and test functions. Computation of the integral for edges belonging to the right, left or bottom boundary are computed in a similar way and omitted here for brevty.

The reaction term, i.e., the last term on the left hand side of \eqref{eq:weak_form_final_disc_coupled} is
\begin{equation}\label{eq:reac_DGDG_coupled}
    \sum_{j=0}^{p_x}\sum_{i=0}^{p_z}q_{m_x,m_z}^{(j,i)}\int_{x_{m_x-1}}^{x_{m_x}}\int_{z_{m_z-1}}^{z_{m_z}}\gamma(x,z)\phi_j^{m_x}(x)\phi_i^{m_z}(z)\phi_{j^\prime}^{m_x}(x)\phi_{i^\prime}^{m_z}(z)\,dxdz
\end{equation}
for $m_x=1,\dots,N_x$, $m_z=1,\dots,N_z$ and
\begin{equation}\label{eq:reac_DGLAG_coupled}
    \sum_{j=0}^{p_x}\sum_{i=0}^{M}q_{m_x,\infty}^{(j,i)}\int_{x_{m_x-1}}^{x_{m_x}}\int_{L_z}^{+\infty}\gamma(x,z)\phi_j^{m_x}(x)\phi_i^{\infty}(z)\phi_{j^\prime}^{m_x}(x)\phi_{i^\prime}^{\infty}(z)\,dxdz.
\end{equation}
If $\gamma$ is constant, then these integrals can be simplified by exploiting orthogonality of the basis functions, otherwise they have to be computed numerically.

Finally, we consider the right hand side of \eqref{eq:weak_form_final_disc_coupled}, which we report here for convenience:
\begin{equation}\label{eq:rhs_coupled}
    \int_{\Omega} f\varphi_h dx\,dz + \sum_{e\in\Gamma_D}\int_e\left(\epsilon\bm{\mu}\nabla{\varphi_h}_{|_e}\cdot\mathbf{n}+\frac{\sigma}{\lvert e\rvert}{\varphi_h}_{|_e}\right)q_D\,ds+\sum_{e\in\Gamma_N}\int_eq_N{\varphi_h}_{|_e}ds.
\end{equation}
The first integral is
\begin{equation}\label{eq:rhs_DGDG_coupled_1st}
    \int_{K_{m_x,m_z}}f(x,z,t)\phi_{j^\prime}^{m_x}(x)\phi_{i^\prime}^{m_z}(z)\,dx\,dz \qquad m_x=1,\dots,N_x \quad m_z=1,\dots,N_z,\infty
\end{equation}
and it is computed using suitable quadrature rules.
The second and third integrals take different forms depending on whether the edge $e$ belongs to the bottom, right or left boundary. If $m_z=1,\dots,N_z$ we have, respectively,
\begin{align}\label{eq:rhs_DGDG_coupled_2nd}
\begin{split}
    &\int_{x_{m_x-1}}^{x_{m_x}}\left(-\mu_z\epsilon\phi_{j^\prime}^{m_x}(x)(\phi_{i^\prime}^\prime)^L+\frac{\sigma}{\Delta x}\phi_{j^\prime}^{m_x}(x)\phi_{i^\prime}^L\right)q_D(x,0,t)\,dx \\
    &\int_{z_{m_z-1}}^{z_{m_z}}\left(\mu_x\epsilon(\phi_{j^\prime}^\prime)^R\phi_{i^\prime}^{m_z}(z)+\frac{\sigma}{\Delta z}\phi_{j^\prime}^R\phi_{i^\prime}^{m_z}(z)\right)q_D(L_x,z,t)\,dz \\
    &\int_{z_{m_z-1}}^{z_{m_z}}\left(-\mu_x\epsilon(\phi_{j^\prime}^\prime)^L\phi_{i^\prime}^{m_z}(z)+\frac{\sigma}{\Delta z}\phi_{j^\prime}^L\phi_{i^\prime}^{m_z}(z)\right)q_D(0,z,t)\,dz
\end{split}
\end{align}
for the second integral and
\begin{align}\label{eq:rhs_DGDG_coupled_3rd}
\begin{split}
    &\int_{x_{m_x}-1}^{x_{m_x}}q_N(x,0,t)\phi_{j^\prime}^{m_x}(x)\phi_{i^\prime}^L\,dx \\
    &\int_{z_{m_z}-1}^{z_{m_z}}q_N(L_x,z,t)\phi_{j^\prime}^R\phi_{i^\prime}^{m_z}(z)\,dz \\
    &\int_{z_{m_z}-1}^{z_{m_z}}q_N(0,z,t)\phi_{j^\prime}^L\phi_{i^\prime}^{m_z}(z)\,dz
\end{split}
\end{align}
for the third. If instead $m_z=\infty$, we only consider the right and the left boundary, for which we obtain
\begin{align}\label{eq:rhs_DGLAG_coupled_2nd}
\begin{split}
    &\int_{L_z}^{+\infty}\left(\mu_x\epsilon(\phi_{j^\prime}^\prime)^R\phi_{i^\prime}^\infty(z)+\sigma\beta\phi_{j^\prime}^R\phi_{i^\prime}^\infty(z)\right)q_D(L_x,z,t)\,dz \\
    &\int_{L_z}^{+\infty}\left(-\mu_x\epsilon(\phi_{j^\prime}^\prime)^L\phi_{i^\prime}^\infty(z)+\sigma\beta\phi_{j^\prime}^L\phi_{i^\prime}^\infty(z)\right)q_D(0,z,t)\,dz
\end{split}
\end{align}
for the second integral and
\begin{align}\label{eq:rhs_DGLAG_coupled_3rd}
\begin{split}
    &\int_{L_z}^{+\infty}q_N(L_x,z,t)\phi_{j^\prime}^R\phi_{i^\prime}^\infty(z)\,dz \\
    &\int_{L_z}^{+\infty}q_N(0,z,t)\phi_{j^\prime}^L\phi_{i^\prime}^\infty(z)\,dz
\end{split}
\end{align}
for the third, respectively.

\section{Local and coupling matrices}
\label{sec: appendix_m}

We report here the expressions for the local matrices in the XDG-Laguerre formulations. For vertical edges $e\in\Gamma^{DG-DG}$, i.e. expression \eqref{eq:diff_internal_DGDG_vertical_coupled}, we have, for $j,j^\prime=0,\dots,p_x$,
\begin{align}
    &m_{j,j^\prime}^{V,11}=-\frac{\mu_x\Delta z}{2}(\phi_j^\prime)^R\phi_{j^\prime}^R+\frac{\mu_x\epsilon\Delta z}{2}\phi_j^R(\phi_{j^\prime}^\prime)^R+\sigma\phi_j^R\phi_{j^\prime}^R, \;
    &m_{j,j^\prime}^{V,12}=\frac{\mu_x\Delta z}{2}(\phi_j^\prime)^R\phi_{j^\prime}^L+\frac{\mu_x\epsilon\Delta z}{2}\phi_j^R(\phi_{j^\prime}^\prime)^L-\sigma\phi_j^R\phi_{j^\prime}^L \\
    &m_{j,j^\prime}^{V,21}=-\frac{\mu_x\Delta z}{2}(\phi_j^\prime)^L\phi_{j^\prime}^R-\frac{\mu_x\epsilon\Delta z}{2}\phi_j^L(\phi_{j^\prime}^\prime)^R-\sigma\phi_j^L\phi_{j^\prime}^R, \;
    &m_{j,j^\prime}^{V,22}=\frac{\mu_x\Delta z}{2}(\phi_j^\prime)^L\phi_{j^\prime}^L-\frac{\mu_x\epsilon\Delta z}{2}\phi_j^L(\phi_{j^\prime}^\prime)^L+\sigma\phi_j^L\phi_{j^\prime}^L.
\end{align}
For horizontal edges $e\in\Gamma^{DG-DG}$, expression \eqref{eq:diff_internal_DGDG_horizontal_coupled}, we have, for $i,i^\prime=0,\dots,p_z$,
\begin{align}
    &m_{i,i^\prime}^{H,11}=-\frac{\mu_z\Delta x}{2}(\phi_i^\prime)^R\phi_{i^\prime}^R+\frac{\mu_z\epsilon\Delta x}{2}\phi_i^R(\phi_{i^\prime}^\prime)^R+\sigma\phi_i^R\phi_{i^\prime}^R, \;
    &m_{i,i^\prime}^{H,12}=\frac{\mu_z\Delta x}{2}(\phi_i^\prime)^R\phi_{i^\prime}^L+\frac{\mu_z\epsilon\Delta x}{2}\phi_i^R(\phi_{i^\prime}^\prime)^L-\sigma\phi_i^R\phi_{i^\prime}^L \\
    &m_{i,i^\prime}^{H,21}=-\frac{\mu_z\Delta x}{2}(\phi_i^\prime)^L\phi_{i^\prime}^R-\frac{\mu_z\epsilon\Delta x}{2}\phi_i^L(\phi_{i^\prime}^\prime)^R-\sigma\phi_i^L\phi_{i^\prime}^R, \;
    &m_{i,i^\prime}^{H,22}=\frac{\mu_z\Delta x}{2}(\phi_i^\prime)^L\phi_{i^\prime}^L-\frac{\mu_z\epsilon\Delta x}{2}\phi_i^L(\phi_{i^\prime}^\prime)^L+\sigma\phi_i^L\phi_{i^\prime}^L
\end{align}
For the vertical edges $e\in\Gamma^{DG-LAG}$ we have: 
where, for $j,j^\prime=0,\dots,p_x$,
\begin{align}
    &m_{j,j^\prime}^{V,11}=-\frac{\mu_x}{2\beta}(\phi_j^\prime)^R\phi_{j^\prime}^R+\frac{\mu_x\epsilon}{2\beta}\phi_j^R(\phi_{j^\prime}^\prime)^R+\sigma\phi_j^R\phi_{j^\prime}^R, \;
    &m_{j,j^\prime}^{V,12}=\frac{\mu_x}{2\beta}(\phi_j^\prime)^R\phi_{j^\prime}^L+\frac{\mu_x\epsilon}{2\beta}\phi_j^R(\phi_{j^\prime}^\prime)^L-\sigma\phi_j^R\phi_{j^\prime}^L \\
    &m_{j,j^\prime}^{V,21}=-\frac{\mu_x}{2\beta}(\phi_j^\prime)^L\phi_{j^\prime}^R-\frac{\mu_x\epsilon}{2\beta}\phi_j^L(\phi_{j^\prime}^\prime)^R-\sigma\phi_j^L\phi_{j^\prime}^R, \;
    &m_{j,j^\prime}^{V,22}=\frac{\mu_x}{2\beta}(\phi_j^\prime)^L\phi_{j^\prime}^L-\frac{\mu_x\epsilon}{2\beta}\phi_j^L(\phi_{j^\prime}^\prime)^L+\sigma\phi_j^L\phi_{j^\prime}^L
\end{align}
where we chose $\lvert e\rvert=1/\beta$. These local matrices are the same as those described earlier in the Appendix; the only difference is that $\Delta z$ has been replaced by $1/\beta$.
Finally, for edges $e\in\Gamma_I^{INTERF}$ the coupling between the discretizations is described by the matrices:
\begin{align}
    &m_{i,i^\prime}^{C,11}=-\frac{\mu_z\Delta x}{2}(\phi_i^\prime)^R\phi_{i^\prime}^R+\frac{\mu_z\epsilon\Delta x}{2}\phi_i^R(\phi_{i^\prime}^\prime)^R+\sigma\phi_i^R\phi_{i^\prime}^R \qquad i=0,\dots,p_z \quad i^\prime=0,\dots,p_z \\
    &m_{i,i^\prime}^{C,12}=\frac{\mu_z\Delta x}{2}(\phi_i^\prime)^R+\frac{\mu_z\epsilon\Delta x}{2}\phi_i^R\left(-\frac{\beta}{2}-\beta i^\prime\right)-\sigma\phi_i^R \qquad i=0,\dots,p_z \quad i^\prime=0,\dots,M \\
    &m_{i,i^\prime}^{C,21}=-\frac{\mu_z\Delta x}{2}\left(-\frac{\beta}{2}-\beta i\right)\phi_{i^\prime}^R-\frac{\mu_z\epsilon\Delta x}{2}(\phi_{i^\prime}^\prime)^R-\sigma\phi_{i^\prime}^R \qquad i=0,\dots,M \quad i^\prime=0,\dots,p_z \\
    &m_{i,i^\prime}^{C,22}=\frac{\mu_z\Delta x}{2}\left(-\frac{\beta}{2}-\beta i\right)-\frac{\mu_z\epsilon\Delta x}{2}\left(-\frac{\beta}{2}-\beta i^\prime\right)+\sigma \qquad i-0,\dots,M \quad i^\prime=0,\dots,M
\end{align}
where we used the properties $\phi_i^\infty(L_z)=1$ and $(\phi_i^\infty)^\prime(L_z)=-\beta(1/2+i)$.

\section{Unknowns and trial functions in flux integral}
\label{sec: appendix_flux_int}

In the sixth integral of \eqref{eq:weak_form_final_disc_coupled}, i.e. the flux integral, expressions for the discrete unknowns and test functions read as follows. Starting with the vertical edges, i.e. expression \eqref{eq:adv_ver_interfaces}, if $e\in\Gamma_I^{DG-DG}$ then $e=\{x_{m_x}\}\times[z_{m_z-1},z_{m_z}]$ and
\begin{align}
    &q_h^a(z)=\sum_{j=0}^{p_x}\sum_{i=0}^{p_z}q_{m_x,m_z}^{(j,i)}\phi_{j}^R\phi_i^{m_z}(z) \qquad q_h^b(z)=\sum_{j=0}^{p_x}\sum_{i=0}^{p_z}q_{m_x+1,m_z}^{(j,i)}\phi_{j}^L\phi_i^{m_z}(z) \\
    &\varphi_h^a(z)=\phi_{j^\prime}^R\phi_{i^\prime}^{m_z}(z) \qquad \varphi_h^b(z)=\phi_{j^\prime}^L\phi_{i^\prime}^{m_z}(z)
\end{align}
while if $e\in\Gamma_I^{DG-LAG}$, then $e=\{x_{m_x}\}\times[L_z,\infty)$ and
\begin{align}
    &q_h^a(z)=\sum_{j=0}^{p_x}\sum_{i=0}^{M}q_{m_x,\infty}^{(j,i)}\phi_{j}^R\phi_i^{\infty}(z) \qquad q_h^b(z)=\sum_{j=0}^{p_x}\sum_{i=0}^{M}q_{m_x+1,\infty}^{(j,i)}\phi_{j}^L\phi_i^{\infty}(z) \\
    &\varphi_h^a(z)=\phi_{j^\prime}^R\phi_{i^\prime}^{\infty}(z) \qquad \varphi_h^b(z)=\phi_{j^\prime}^L\phi_{i^\prime}^{\infty}(z).
\end{align}
In both cases $F_1^a(q_h^a)=F_1(q_h^a,x_{m_x}^-,z)$ and $F_1^b(q_h^b)=F_1(q_h^b,x_{m_x}^+,z)$.

For horizontal edges, expression \eqref{eq:adv_hor_interfaces},
If $m_z<N_z$, then
\begin{align}
    &q_h^a(x)=\sum_{j=0}^{p_x}\sum_{i=0}^{p_z}q_{m_x,m_z}^{(j,i)}\phi_{j}^{m_x}(x)\phi_i^{R} \qquad q_h^b(x)=\sum_{j=0}^{p_x}\sum_{i=0}^{p_z}q_{m_x,m_z+1}^{(j,i)}\phi_{j}^{m_x}(x)\phi_i^{L} \\
    &\varphi_h^a(x)=\phi_{j^\prime}^{m_x}(x)\phi_{i^\prime}^{R} \qquad \varphi_h^b(x)=\phi_{j^\prime}^{m_x}(x)\phi_{i^\prime}^{L}
\end{align} 
while if $m_z=N_z$ we obtain
\begin{align}
    &q_h^a(x)=\sum_{j=0}^{p_x}\sum_{i=0}^{p_z}q_{m_x,m_z}^{(j,i)}\phi_{j}^{m_x}(x)\phi_i^{R} \qquad q_h^b(x)=\sum_{j=0}^{p_x}\sum_{i=0}^{M}q_{m_x,\infty}^{(j,i)}\phi_{j}^{m_x}(x) \\
    &\varphi_h^a(x)=\phi_{j^\prime}^{m_x}(x)\phi_{i^\prime}^{R} \qquad \varphi_h^b(x)=\phi_{j^\prime}^{m_x}(x)
\end{align}
where we used $\phi_i^\infty(L_z^+)=0$ for all $i$.

\section{Additional numerical results for convergence tests and coupling validation}\label{sec:tables}

In this section we report detailed results related to the convergence tests and the coupling validation tests, see Sections  \ref{subsec:conv_tests} \ref{subsec:coupling_valid} in the main text.

\begin{table}[ht]
\centering
\begin{tabular}[]{ccc}\toprule
$M$ & $\mathcal{E}_2^{abs}$ & $\mathcal{E}_\infty^{abs}$\\ \midrule
5 & 3.13e-04 & 2.17e-03\\
10 & 1.05e-04 & 1.68e-04 \\
20 & 1.46e-06 & 9.79e-05 \\
30 & 2.61e-08 & 1.31e-07 \\
35 & 2.25e-09 & 1.68e-08\\\bottomrule
\end{tabular}\caption{Convergence with respect to $M$ of the 2D DG-Laguerre scheme, linear advection-diffusion equation on a semi-infinite strip. $L_x=1\,m$, $N_x=1000$, $p_x=3$, $\beta=5$, $T=5\times 10^{-3}\,s$, $N_t=10$, $\mu_x=0.05\,m^2/s$, $\mu_z=0.01\,m^2/s$, $u_x=1\,m/s$, $u_z=2\,m/s$.}\label{tab:M_N300}
\end{table}

\begin{table}[ht]
\centering
\begin{tabular}{ccccc}\toprule
$N_x$,$p_x=1$ & $\mathcal{E}_2^{abs}$ & $r_2$ & $\mathcal{E}_\infty^{abs}$ & $r_\infty$ \\ \midrule
40 & 1.86e-03 & & 5.07e-03 &  \\
80 & 5.31e-04 & 1.81 & 1.46e-03 & 1.79 \\
160 & 1.43e-04 & 1.90 & 3.89e-04 & 1.91 \\
320 & 3.70e-05 & 1.95 & 1.01e-04 & 1.95\\ \bottomrule
\end{tabular}
\qquad
\begin{tabular}{ccccc}\toprule
$N_x$,$p_x=2$ & $\mathcal{E}_2^{abs}$ & $r_2$ & $\mathcal{E}_\infty^{abs}$ & $r_\infty$ \\ \midrule
40 & 9.66e-04 & & 2.89e-03 &  \\
80 & 2.99e-04 & 1.69 & 8.80e-04 & 1.71 \\
160 & 8.18e-05 & 1.87 & 2.41e-04 & 1.87 \\
320 & 2.12e-05 & 1.94 & 6.26e-05 & 1.95\\ \bottomrule
\end{tabular}

\bigskip

\begin{tabular}{ccccc}\toprule
$N_x$,$p_x=3$ & $\mathcal{E}_2^{abs}$ & $r_2$ & $\mathcal{E}_\infty^{abs}$ & $r_\infty$ \\ \midrule
50 & 2.03e-06 & & 7.45e-06 & \\
75 & 3.93e-07 & 4.05 & 1.42e-06 & 4.09 \\
100 & 1.24e-07 & 4.02 & 4.40e-07 & 4.07 \\
125 & 5.05e-08 & 4.01 & 1.78e-07 & 4.06 \\
150 & 2.43e-08 & 4.00 & 8.50e-08 & 4.05 \\
175 & 1.31e-08 & 4.00 & 4.56e-08 & 4.03\\\bottomrule
\end{tabular}\caption{Absolute $L^2$ and $L^\infty$ errors and convergence rates, $x$-convergence, $p_x=1$, $p_x=2$ and $p_x=3$.}\label{tab:x_err_eps1}
\end{table}

\begin{table}[ht]
\centering
\begin{tabular}[]{cccccc}\toprule
$N_x$ & $N_t$ & $\mathcal{E}_2^{abs}$ & $r_2$ & $\mathcal{E}_\infty^{abs}$ & $r_\infty$\\ \midrule
50 & 100 & 9.80e-04 & & 2.91e-03 & \\
100 & 200 & 2.74e-04 & 1.84 & 7.99e-04 & 1.86 \\
150 & 300 & 1.26e-04 & 1.91 & 3.65e-04 & 1.93 \\
200 & 400 & 7.25e-05 & 1.93 & 2.09e-04 & 1.94 \\
250 & 500 & 4.70e-05 & 1.95 & 1.35e-04 & 1.96\\\bottomrule
\end{tabular}\caption{Space-time convergence of the 2D DG-Laguerre scheme, linear advection-diffusion equation on a semi-infinite strip. $L_x=1\,m$, $p_x=1$, $\beta=6$, $M=60$, $T=1\,s$, $\mu_x=0.05\,m^2/s$, $\mu_z=0.01\,m^2/s$, $u_x=1\,m/s$, $u_z=2\,m/s$.}\label{tab:x_t_conv_eps1}
\end{table}

\begin{table}[hb]
\centering
    \begin{tabular}{cccccc}\toprule
        $M$ & $\beta$ & $\sigma_x=\sigma_z$ & $\mathcal{E}_2^{rel}$ & $\mathcal{E}_{\infty}^{rel}$ \\ \midrule
        10 & 4 & 0.5 & 1.46e-02 & 3.73e-02 \\
        & & 1 & 2.87e-03 & 7.53e-03 \\
        & & 2 & 1.28e-04 & 4.22e-04 \\ \midrule
        40 & 6 & 0.5 & 2.18e-04 & 4.89e-04 \\
        & & 1 & 6.44e-05 & 1.78e-04 \\
        & & 2 & 1.39e-06 & 4.74e-06\\\bottomrule
    \end{tabular} \quad \begin{tabular}{cccc}
    \toprule
    $M$ & $\mathcal{E}_2^{rel}$ & $\mathcal{E}_\infty^{rel}$ \\ \midrule
    10 & 1.17e-02 & 7.10e-02 \\
    20 & 1.43e-03 & 1.01e-02 \\
    40 & 1.85e-04 & 1.55e-03 \\
    80 & 9.54e-05 & 4.60e-04\\\bottomrule
    \end{tabular}\caption{Left: 2D Advection-diffusion equation, Gaussian initial data. Relative errors at final time  $T=4\,s$ in $[0,10\,m]\times[0,10\,m]$ of the XDG discretization with respect to a single-domain DG discretization on $[0,10\,m]\times[0,20\,m]$. $L_x=L_z=10\,m$, $p_x=p_z=1$, $N_x=50$, $N_z=500$, $N_t=200$, $\mu_x=\mu_z=0.1\,m^2/s$, $u_x=0.5\,m/s$, $u_z=1\,m/s$. Right: Burgers' equation. Relative errors in $[0,10\,m]\times[0,10\,m]$ of the XDG discretization with respect to a single-domain DG discretization on $[0,10\,m]\times[0,20\,m]$. $L_x=L_z=10\,m$, $p_x=p_z=1$, $N_x=50$, $N_z=80$, $\beta=10$, $T=5\,s$, $N_t=500$, $\mu_x=\mu_z=0.05\,m^2/s$.}\label{tab:adv_diff_coup_valid}
\end{table}

\newpage\clearpage
\bibliography{dg_laguerre_2d}   

\end{document}